\documentclass[11pt]{article}  
\usepackage{amsmath}
\usepackage{amssymb}
\usepackage{bm}
\usepackage{multirow}
\usepackage{setspace}
\usepackage{comment}
\usepackage{graphicx,epstopdf}
\usepackage[margin=1 in]{geometry}
\usepackage[font=small, labelfont=bf]{caption}
\usepackage{url}

\newcommand{\beqnal}{\begin{equation}\begin{aligned}}
\newcommand{\eeqnal}{\end{aligned} \end{equation}}
\newcommand{\beqn}{\begin{eqnarray}}
\newcommand{\bdm}{\begin{displaymath}}
\newcommand{\edm}{\end{displaymath}}
\newcommand{\eeqn}{\end{eqnarray}}
\newcommand{\be}{\begin{equation}}
\newcommand{\ee}{\end{equation}}
\newcommand{\ba}{\begin{array}}
\newcommand{\ea}{\end{array}}

\usepackage[super]{natbib}
\usepackage[small,compact]{titlesec}

\makeatletter
\renewcommand\@biblabel[1]{#1.}
\makeatother

\def\R{{\rm I\kern-.19em R}}
\def\M{{\rm I\kern-.1567em M}}
\def\C{{\rm I\kern-.55em C}}

\def\Chi{{\raise.4ex\hbox{\large$\chi$}}}

\author{Calvin Tsay and Michael Baldea\thanks{corresponding author: mbaldea@che.utexas.edu}\\
	McKetta Department of Chemical Engineering, The University of Texas at Austin\\
	Austin, TX 78712}

\title{Integrating Production Scheduling and Process Control using Latent Variable Dynamic Models}

\titleformat{\section}
{\normalfont\bfseries\normalsize}
{\thesection}{1em}{}
\titleformat{\subsection}
{\normalfont\mdseries\bfseries\normalsize}
{\thesubsection}{1em}{}
\titlespacing{\subsection}
{2pc}{0.35ex plus .1ex minus .2ex}{0.3ex minus .1ex}

\titleformat{\subsubsection}
{\normalfont\itshape\mdseries\normalsize}
{\thesubsubsection}{1em}{}
\titlespacing{\subsubsection}
{4pc}{0.5ex plus .1ex minus .2ex}{0.3ex minus .1ex}

\begin{document}
\singlespacing
\maketitle

\begin{abstract}
With the increasing need for dynamic decisions in fast-changing markets, the integration of scheduling and control is an important consideration in chemical processes. Nevertheless, computing optimal production schedules using dynamic process models remains challenging, due to model nonlinearity and high-dimensionality. In this paper, we observe that the intrinsic dimensionality of process dynamics (as relevant to scheduling) is often much lower than the number of model state and/or algebraic variables. We introduce a data mining approach to ``learn'' closed-loop process dynamics on a low-dimensional, latent manifold. The manifold dimensionality is selected based on a tradeoff between model accuracy and complexity. After projecting process data, system identification and optimal scheduling calculations can be performed in the low-dimensional, latent-variable space. We apply these concepts to schedule an air separation unit under time-varying electricity prices. We show that our approach reduces the computational effort required, while offering more detailed dynamic information compared to previous related works.
\end{abstract}


\section{Introduction}
The fast-changing conditions present in modern markets introduce many challenges and opportunities for improvement in chemical process operations. For example, partly due to increased adoption of intermittent renewable energy sources, real-time electricity prices can fluctuate by several orders of magnitude during a 24-hour period. This provides a strong incentive for demand response, or intentional modification by an electricity user of its power consumption (``load shifting'') over time in order to exploit time-dependent electricity prices \cite{zhang2015}. In such circumstances, production scheduling must consider production changes over a sequence of relatively short time slots (intervals) to maximize profits. To ensure the resulting production schedules are feasible when implemented in the physical process, many research works and applications seek to integrate decision-making across different levels and time-scales\cite{dowling2017, otashu2018, schafer2018}, with a focus on accounting for process agility (expressed in terms of dynamic characteristics and/or control performance) in scheduling calculations. Thus, the \textit{integration of scheduling and control} for chemical processes has emerged as an important research area \cite{engell2012, baldea2014, daoutidis2018, dias2019}.

Production scheduling refers to the determination of production sequences, product grades, batch sizes, unit assignments, and/or task timing that maximize profits (or minimize cost). Scheduling decisions are typically made over a time horizon spanning several hours to several days. Most of the conventional methods for computing optimal schedules rely on the assumption that the process is at a steady state before each change in production targets, and that it will again reach a new steady state shortly thereafter. However, this assumption may not be valid when scheduling decisions are made over shortened time intervals, such as those required by the aforementioned, fast-changing market conditions \cite{baldea2014}. Once the optimal schedule is determined, the control layer of a process seeks to track the setpoints/targets determined by the scheduling layer, while satisfying process and product constraints. To this end, optimization-based controllers (notably model predictive control--MPC) have enjoyed widespread acceptance in the chemical industry \cite{qin2003}. MPC determines optimal control moves using a dynamic model to predict the plant response (over a prediction horizon typically much shorter than the horizon of scheduling calculations) to changes in the manipulated variables.

In an early effort to integrate process scheduling and control, Flores-Tlacuahuac and Grossmann \cite{flores2006} explicitly included the dynamic process model and controller in the scheduling problem, resulting in a large simultaneous dynamic optimization problem. Zhuge and Ierapetritou \cite{zhuge2012} later implemented this discretized-time approach in a closed-loop strategy to mitigate the effect of disturbances. Beal et al. \cite{beal2018} extended this concept to account for time-dependent parameters and constraints, and demonstrated the economic benefits of integrated scheduling control in both open-loop and closed-loop implementations \cite{beal2017}. Koller et al. \cite{koller2018} considered embedding PI controllers into scheduling calculations, accounting for stochastic disturbances and uncertainties using a sample-based, back-off method.

In general, embedding a dynamic process model in scheduling calculations tends to increase computational cost, and many optimization techniques have been proposed to facilitate dealing with the integrated problem. Nystr\"{o}m et al. \cite{nystrom2005} reduced the computational complexity of solving the integrated problem by decomposing it into a scheduling master problem and a control sub-problem, and iterating between the two. Nie et al. \cite{nie2014} took a similar approach, proposing a generalized Benders decomposition algorithm, where the scheduling decisions comprise the master problem and the dynamic process optimization comprises the primal problem.
Simkoff and Baldea \cite{simkoff2019} directly incorporated the KKT optimality conditions of a linear MPC system using complementary constraints to provide an exact representation of closed-loop dynamics in the scheduling layer.

Several works \cite{burnak2018, zhuge2014} employed multiparametric model predictive control (mpMPC) in an optimal scheduling framework. Charitopoulos et al. \cite{charitopoulos2018} examined the closed-loop implementation of an mpMPC approach for integrated scheduling and control that can handle dynamic disturbances. mpMPC approaches rely on generating explicit forms of the optimal control laws offline, and are thus computationally efficient when implemented online. However, the complexity of the offline problem grows exponentially with the size of the system (model) under consideration and with the dimension of its input and output vectors.

For large-scale process models, the computational cost of optimal scheduling calculations using dynamic models can be diminished by model reduction, using, e.g.,  collocation- or compartmentalization-based methods \cite{cao2016, schafer2019}. Similarly, Du et al. \cite{du2015} reduced the order of closed-loop process dynamics with input-output feedback linearizing controllers, designed to impose a well-defined linear closed-loop response that could be easily represented in scheduling calculations. Baldea et al. \cite{baldea2015} extended this approach for systems under MPC by modifying the MPC constraint set, such that the controller imposes the desired linear closed-loop behavior. While the above model reduction approaches rely on first-principles dynamic models, recent works \cite{pattison2016, tsay2019} have proposed deriving reduced-order approximations of process dynamics from recorded historical operating data via system identification. However, in this latter context, selecting the structure and scope of a scheduling-relevant, reduced-order model required significant engineering expertise. Thus, the development of a systematic  framework and corresponding automated workflow for performing scheduling-oriented modeling and/or model reduction of the dynamics of chemical processes remains an open problem.

Motivated by the above, in this work, we present a framework for systematically deriving low-dimensional, scheduling-relevant dynamic models using a latent-variable representation. We first analyze the intrinsic (approximate) low-dimensionality of closed-loop process dynamics relevant to scheduling. This observation motivates learning the underlying \textit{latent manifold} that describes the process behavior in its intrinsic dimension, and using the transformation to improve the computational performance of scheduling calculations. The novelty of this contribution consists of the following:
\begin{itemize}
	\item A conceptual approach for identifying a low-dimensional manifold underlying closed-loop process behavior using data-mining techniques. Specifically, we show that autoencoders provide a simple method to learn nonlinear relationships among process variables.
	\item A framework for system identification in the space of the learned latent variables. In contrast to previous works \cite{pattison2016, tsay2018, tsay2019} that rely on physical insight to select a subset of ``scheduling-relevant'' variables from process state and output vectors, the present framework automates the dimensionality-reduction process, eliminating empiricism from this step.
	\item A methodology for integrated scheduling and control using the above low-dimensional representation of closed-loop process dynamics. Owing to its low intrinsic dimensionality, the resulting problem has low computational complexity, while retaining sufficiently rich dynamic information.
\end{itemize}
We apply the proposed approach to derive data-driven, latent-variable models of an air separation unit (ASU) and compute optimal schedules for a demand-response scenario where electricity prices fluctuate on an hourly basis. We show that the results compare favorably against previous works \cite{pattison2016, dias2018} that identify low-order dynamic models of an empirically defined subset of process variables. Finally, we investigate tuning the dimensionality of the reduced-order representation (i.e., selecting the number of latent variables) to manage the tradeoff between optimization problem size and model accuracy.

\section{Background}
\subsection{Integrated Scheduling and Control Problem}
The integrated scheduling and control problem aims to derive dynamically feasible production schedules by including the \textit{closed-loop} dynamics of a process in the scheduling problem. The deterministic, continuous-time integrated scheduling and control problem can be stated generally as \cite{du2015}:
\begin{eqnarray}
\underset{\mathbf{y}_{sp}(t)}{\text{max}} &\int_{t=0}^{t=t_f} P(\mathbf{y}, t) dt \label{eq:scheduling1}\\
s.t. & \dot{\mathbf{x}} = f(\mathbf{x}) + G(\mathbf{x})\mathbf{u} \\
&\mathbf{y} = h(\mathbf{x}) \\
&\mathbf{u} = K(\mathbf{y}_{sp}-\mathbf{y}) \\
&l(\mathbf{x}, \mathbf{y}, \mathbf{u}, t) = 0 \\
&\Bigg[ \begin{array}{c} \mathbf{x}^L \\ \mathbf{y}^L \\ \mathbf{u}^L \end{array}\Bigg]
\leq \Bigg[ \begin{array}{c} \mathbf{x} \\ \mathbf{y} \\ \mathbf{u} \end{array}\Bigg]
\leq \Bigg[ \begin{array}{c} \mathbf{x}^U \\ \mathbf{y}^U \\ \mathbf{u}^U \end{array}\Bigg] \label{eq:scheduling2}
\end{eqnarray}
where $\mathbf{y}_{sp}(t)$ is a time-varying vector of production targets and/or other setpoints to be supplied to the control system, and $K$ represents the process control policy. The vector $\mathbf{x} \subset \mathbb{R}^n$ denotes the process state variables, $\mathbf{y} \subset \mathbb{R}^m$ are the output variables, $\mathbf{u} \subset \mathbb{R}^u$ are the input variables, $f$ and $h$ are appropriately defined vector fields, and $G$ is of appropriate dimensions. The process variables $\mathbf{x}, \mathbf{y}, \mathbf{u}$ may be subject to constraints, given by \eqref{eq:scheduling2}. The economic objective function $P(\mathbf{y},t)$ typically includes the revenue from selling product, and the process operating costs. $l(\cdot)$ includes storage and demand constraints that ensure (i) customer demand can be met at all times, (ii) the amount of product stored does not deplete/exceed the physical capacity of the storage system, and (iii) artificial economic gains are not realized by depleting material hold-up present at $t=0$. Note that $l(\cdot)$ may include both path ($0 \leq t \leq t_f$) and endpoint ($t=t_f$) constraints.

\subsection{Scale-Bridging Models}
Theoretical developments \cite{baldea2015, baldea2012} have shown that the closed-loop, input-output dynamic behavior of process systems (i.e., the response of the process to changes in $\mathbf{y}_{sp}$) may be quite slow in comparison to the evolution of states of the individual process units. The former, input-output behavior often evolves over time scales relevant to process scheduling calculations, particularly in the context of the fast-changing market conditions. Moreover, these results \cite{baldea2015, baldea2012} suggest that the input-output dynamic behavior can be described, or at least usefully approximated, using a low-order model of dimension much smaller than that of the state variable vector $\mathbf{x}$. On this basis, our previous works proposed representing the process dynamics in scheduling calculations using time scale-bridging models (SBMs), which are a low-order representation of the closed-loop dynamics of the process. Broadly speaking, low-order SBMs can be derived using either model reduction or system identification, and the two techniques are briefly reviewed below:
\begin{itemize}
\item \textbf{Model reduction} refers to systematic derivation of a low-order model from a detailed (high-dimensional and likely first-principles) dynamic process model. Many methods have been proposed for systematic model reduction, such as asymptotic analyses \cite{baldea2012, kumar1998} and null-space projection \cite{yu2015}. These methods reduce the number of states, typically in systems that exhibit multiple time scale dynamics, resulting in a lower-dimensional differential (algebraic) equation system. These modeling techniques retain physically meaningful states, which can be a useful feature for controller design and monitoring strategies \cite{baldea2012, yu2016}. Alternative techniques such as empirical Gramians \cite{hahn2002} and proper orthogonal decomposition \cite{willcox2002, armaou2002} can also be used, although these methods can result in states that are not physically meaningful.

\item \textbf{System identification} refers to ``learning'' a dynamic model from process operating data. In this approach, a full-order, dynamic process model is not required. Rather, a generic model structure is assumed, and its parameters are computed based on transient data recorded from the process. The data can be generated in the course of the operation of a physical process (either routine or following a deliberate pattern of system identification experiments) or from simulations of a detailed mathematical model. Learning (and reducing the dimensionality of) dynamics from data is a highly active area of research \cite{klus2018}. We direct the interested reader to, e.g., the book by Zhu \cite{zhu2001} for an overview of the multitude of available system identification techniques and their application to chemical processes.
\end{itemize}
Focusing on system identification methods, Pattison et al. \cite{pattison2016} suggested that the dimensionality (and consequently the computational complexity) of the scheduling problem \eqref{eq:scheduling1}--\eqref{eq:scheduling2} could be considerably reduced by restricting modeling efforts to ``scheduling-relevant variables.'' In this approach, a subset $\mathbf{z} \subseteq [\mathbf{x}, \mathbf{y}, \mathbf{u}]$ is defined, that includes the input and output variables ($\mathbf u$ and $\mathbf y$) that affect $P(\cdot)$ and $l(\cdot)$. Process state and output variables ($\mathbf x$ and $\mathbf y$) whose constraints are active during steady-state operation or during process transitions are intuitively assumed to limit the process dynamic agility and are also included in $\mathbf{z}$ to ensure that the resulting schedules do not violate any constraints (i.e., guarantee the dynamic feasibility of a schedule). In essence, the dimensionality reduction occurs through a heuristic and expertise-intensive selection of the scheduling-relevant variables $\hat{\mathbf{x}} \subseteq \hat{\mathbf{x}}$, $\hat{\mathbf{y}} \subseteq \hat{\mathbf{y}}$, and $\hat{\mathbf{u}} \subseteq \hat{\mathbf{u}}$. The resulting scheduling optimization problem is similar to \eqref{eq:scheduling1}--\eqref{eq:scheduling2}, but can contain significantly fewer variables and constraints:
\begin{eqnarray}
\underset{\mathbf{y_{sp}}}{\text{max}} &\int_{t=0}^{t=t_f} P(\mathbf{z}, t) dt \label{eq:sbm1}\\
s.t. &\dot{\mathbf{z}} = \Bigg[\begin{array}{c} \dot{\hat{\mathbf{y}}} \\ \dot{\hat{\mathbf{x}}} \\ \dot{\hat{\mathbf{u}}} \end{array} \Bigg] = \hat{f}(\mathbf{z}, \mathbf{y}_{sp}) \\
&l(\mathbf{z},t) = 0 \\
&\mathbf{z}^{L} \leq \mathbf{z} \leq \mathbf{z}^{U} \label{eq:sbm2}
\end{eqnarray}
where $l(\mathbf{z},t)$ may include both path and endpoint storage/demand constraints. The reduced-space dynamic model $\hat{f}(\mathbf{z}, \mathbf{y}_{sp})$ represents the scale-bridging model (SBM), which approximates the closed-loop, input-output relationships between the process setpoints and the scheduling-relevant variables and replaces the process model present in \eqref{eq:scheduling1}--\eqref{eq:scheduling2}.

\subsection{Manifold Learning}
While limiting the scheduling problem to only consider the dynamics of ``scheduling-relevant variables'' can be an effective form of dimensionality reduction, the selection of these variables relies on manual effort, technical insight, and human expertise. On the other hand, chemical processes typically have many sensors that record process variables at frequencies in the order of minutes, generating ``big data'' sets that can be exploited to understand the underlying system behavior. There exist many approaches for learning low-dimensional representations of a dynamical system from recorded data. In this context, \textit{manifold learning} refers to identifying a low-dimensional manifold on which higher-dimensional data points intrinsically lie. The learned manifold represents a subspace of the full-dimensional variable space that explains (most of) the variation observed in the data set. Observations of the original system can be transformed to (projected on) a smaller set of \textit{latent variables} that parameterize the manifold.

A broad class of unsupervised machine learning algorithms can be applied to the task of manifold learning. Pearson \cite{pearson1901} introduced principal component analysis (PCA) in 1901, and the technique is now a widely accepted dimensionality-reduction technique. PCA consists of finding a linear coordinate transformation whereby the data are projected on a new set of latent variables. The coordinate transform is constructed such that the amount of variance captured by each successive latent variable, or principal component, is maximized. Latent variables based on linear combinations of the original variables, such as those from PCA, are commonly used in the process industries for process monitoring and troubleshooting \cite{macgregor2005}. They have also found applications in process control, where they can be employed to reduce the dimension of the controlled variable space. For example, latent variables can replace the original process controlled variables to simplify controller calculations \cite{flores2005, lauri2010}.

While PCA is limited to finding linear mappings, a number of nonlinear manifold learning algorithms have been presented. A simple nonlinear extension of PCA is kernel PCA, where a nonlinear kernel is first applied, and PCA is performed in the processed feature space \cite{scholkopf1997}. Several researchers have studied the relationships between PCA and a particular class of artificial neural network known as autoencoders. Of particular note, Sanger \cite{sanger1989} showed that linear autoencoders correspond exactly to PCA, while Kramer \cite{kramer1991} proposed nonlinear autoencoders as a form of generic, nonlinear PCA. Many other nonlinear manifold learning techniques have been since proposed, including diffusion maps \cite{coifman2006}, Laplacian eigenmaps \cite{belkin2003}, locally linear embeddings \cite{roweis2000}, stochastic neighbor embeddings (SNE) \cite{hinton2003}, and Isomap \cite{tenenbaum2000}. For further details, the reader is referred to the book by Lee and Verleysen \cite{lee2007} and the review by Van Der Maaten et al. \cite{van2009}.

\section{Scheduling with Learned Latent Variables}
\subsection{Concept}
We observe now that the \textit{intrinsic dimensionality} (the number of independent variables underlying the significant nonrandom variations in the observations) of the closed-loop behavior of a chemical process can be much lower than the apparent \textit{extrinsic dimensionality} ($n+m+u$). In particular, the dimensionality of the process input $\mathbf{u}$, output $\mathbf{y}$, and state $\mathbf{x}$ variables can be parameterized by the process state variables $\mathbf{x}$:
\begin{eqnarray}
\mathbf{x}^* \equiv \Bigg[ \begin{array}{c} \mathbf{x} \\ \mathbf{y} \\ \mathbf{u} \end{array}\Bigg]
= \Bigg[ \begin{array}{c} \mathbf{x} \\ h(\mathbf{x}) \\ K(\mathbf{y}_{sp}-h(\mathbf{x})) \end{array}\Bigg]
\end{eqnarray}
where we define the augmented process state variable vector $\mathbf{x}^*$ as $[\mathbf{x}, \mathbf{y}, \mathbf{u}]$. Therefore, the mapping relating the setpoints/targets set by the scheduling layer, $\mathbf{y}_{sp}$, to $\mathbf{x}^*$ has an intrinsic dimensionality equal to $\text{dim}_i(\mathbf{x})$. We use $\text{dim}(\cdot)$ to denote extrinsic dimensionality and $\text{dim}_i(\cdot)$ to denote the intrinsic dimensionality, as defined above. Note that this assumes $\mathbf{u} = K(\mathbf{y}_{sp}-h(\mathbf{x}))$ can be evaluated directly, i.e., an explicit control law exists. In the case of an implicit/optimization-based controller, an explicit relationship may still exist, or an approximation may be possible \cite{pistikopoulos2009}. We shall examine the case of dimensionality reduction for a process operating under an optimization-based controller in the study presented later in the paper.

For the particular case of model predictive control (MPC), Lovelett et al. \cite{lovelett2018} observed that the optimal control policy $\mathbf{u}(t)$ may have a significantly lower \textit{intrinsic} dimensionality than its \textit{extrinsic} dimensionality--an observation aligned with the aforementioned findings concerning the process dynamics. The extrinsic dimensionality of $\mathbf{u}(t)$ is equal to $N \times \text{dim}(\mathbf{u})$, where $N$ is the number of computed steps, while the intrinsic dimensionality is limited by $\text{dim}_i(\mathbf{u}) \leq \text{dim}(\mathbf{y}) + \text{dim}(\mathbf{y}_{sp}(t))$. Here, $\mathbf{y}$ represents the vector of current values of the output variables and $\mathbf{y}_{sp}(t))$ represents the reference trajectory. Note that $\text{dim}(\mathbf{y}_{sp}(t))$ is dependent on the control vector parameterization of $\mathbf{y}_{sp}(t)$ with respect to time. The intrinsic dimensionality of the control policy $\text{dim}_i(\mathbf{u})$ is often lower than this limit. This may occur for several reasons, including the above case where the dynamics of the system itself lie on a low-dimensional, slow manifold \cite{baldea2012}, or where the state-space realization is non-minimal order (containing redundant information \cite{lovelett2018}).

In this work, we propose a new learning-based approach for low-order SBM generation, whereby we find a latent manifold mapping of the augmented process state variable vector $\mathbf{x}^*$. We seek an invertible mapping $\mathbf{x}^* \leftrightarrow \bm{\phi}$, with $\bm{\phi} \in \mathbb{R}^p$, denoted as $\bm{\phi} = c(\mathbf{x}^*)$ and $\mathbf{x}^* = c^{inv}(\bm{\phi})$. Furthermore we desire to identify the mapping $c: \mathbb{R}^n \times \mathbb{R}^m \times \mathbb{R}^u \rightarrow \mathbb{R}^p$ such that $p << n+m+u$.  We note that such a mapping always exists for $p \leq n+m+u$, since a trivial exact mapping is possible at $dim(\bm{\phi}) = dim(\mathbf{x}^*)$. Once a mapping is identified, the dynamics of the \textit{latent variables} $\bm{\phi}$ can be embedded in the scheduling problem. The resulting scheduling problem has a low intrinsic dimensionality, with dynamics evolving only in the (low-dimensional) latent-variable space $\mathbb{R}^{p}$:
\begin{eqnarray}
\underset{\mathbf{y}_{sp}}{\text{max}} &\int_{t=0}^{t=t_f} P(\mathbf{y}, t) dt \label{eq:latentsched1} \\
s.t. &\dot{\bm{\phi}} = f^\phi(\bm{\phi}, \mathbf{y}_{sp}) \\
& \mathbf{x}^* \equiv \Bigg[ \begin{array}{c} \mathbf{x} \\ \mathbf{y} \\ \mathbf{u} \end{array}\Bigg] = c^{inv}(\bm{\phi}) \\
& l(\mathbf{x}^*,t) = 0 \\
&\mathbf{x}^{*L} \leq \mathbf{x}^* \leq \mathbf{x}^{*U} \label{eq:latentsched2}
\end{eqnarray}
Assuming that the mapping $c: \mathbf{x} \rightarrow \bm{\phi}$ and the inverse mapping $c^{inv}: \bm{\phi} \rightarrow \mathbf{x}$ exist, the dimensionality of the dynamic constraint(s) is now $p = dim(\bm{\phi})$. If the mappings $c(\cdot)$, $c^{inv}(\cdot)$ are exact, and the dynamics of the latent variables are represented accurately by $f^\phi(\cdot)$, then \eqref{eq:latentsched1}--\eqref{eq:latentsched2} is identical to the original scheduling problem \eqref{eq:scheduling1}--\eqref{eq:scheduling2}. Note that for $c(\cdot)$, $c^{inv}(\cdot)$ to be exact, or equivalently, $c^{inv}(c(\mathbf{x}^*)) = \mathbf{x}^*$, $\mathbf{x}^*$ can only contain $p$ independent variables. The remaining variables must be (nonlinearly) correlated. In practical situations, process variables that feature path constraints of the type in \eqref{eq:latentsched2} may only be a subset of the full vector of process variables, and manifold learning can be carried out in a space of already lower dimension. Nevertheless, recent work \cite{obermeier2019, wiebe2018} has highlighted tradeoffs between dynamic production schedules and equipment fatigue, suggesting that some variables without explicit constraints may still be relevant in the scheduling layer and should be included in $\mathbf{x}^*$.

If the dynamics of $\mathbf{x}^*$ present a low-dimensional manifold only in a limit case (e.g., when the process dynamics are in a singularly perturbed form), the low-dimensional dynamics only approximate the true system. Specifically,  $c^{inv}(c(\mathbf{x}^*)) \approx \mathbf{x}^*$, and the mappings $c(\cdot)$, $c^{inv}(\cdot)$ are inexact. In this case, some information is lost by ``collapsing'' the dynamics of $\mathbf{x}^*$ onto a reduced dimension, resulting in an approximation we denote as ${\mathbf{x}^*}'=c^{inv}(c(\mathbf{x}^*))$. The accuracy of the approximation, in terms of $\parallel\mathbf{x}^* - {\mathbf{x}^{*}}'\parallel$, can be improved by increasing $p$ until the original system is fully recovered at $p = n+m+u$ ($p$ can be smaller if some variables are correlated). In other words, the dimension of the latent manifold $p$ can be used as a parameter for adjusting the accuracy of the reduced-order representation of the closed-loop dynamics.

\textbf{Remark 1.} \textit{Reducing the number of dynamic variables in the original problem \eqref{eq:sbm1}--\eqref{eq:sbm2} to the lower-dimensional problem \eqref{eq:latentsched1}--\eqref{eq:latentsched2} may be beneficial for both sequential \cite{tsay2019} and simultaneous \cite{kelley2018} dynamic optimization approaches. In this work, we focus on sequential approaches, where the Jacobian size for computing implicit time integration steps is reduced by limiting the number of dynamic variables to $p$. We expect the benefits to also extend to simultaneous approaches, where the number of differential state variables treated by the optimization problem is still reduced. However, note that the explicit dimensionality of $dim(\mathbf{x}^*)$ may remain larger than in \eqref{eq:sbm1}--\eqref{eq:sbm2}. We refer the interested reader to \cite{vassiliadis1994, vassiliadis1994b} for an overview of sequential techniques for dynamic optimization and \cite{biegler2007} for information on simultaneous strategies.
}

\subsection{Latent Variable Scheduling Framework}
The proposed approach for latent variable scheduling comprises the following steps:
\begin{enumerate}
	\item Obtain historical process operating data representative of typical production schedules
	\item Learn latent variable mappings $c: \mathbf{x}^* \rightarrow \bm{\phi}$ and $c^{inv}: \bm{\phi} \rightarrow {\mathbf{x}^*}'$
	\item Transform historical data $\mathbf{x}^*(t)$ using $c$ to produce $\bm{\phi}(t)$
	\item Determine model form and fit a dynamic model $\dot{\bm{\phi}} = f^\phi(\bm{\phi},\mathbf{y}_{sp})$ to the latent variables using the transformed data set
	\item Solve the low-dimensional scheduling problem \eqref{eq:latentsched1}--\eqref{eq:latentsched2} with path constraints on ${\mathbf{x}^*}'=c^{inv}(\bm{\phi})$
\end{enumerate}

\subsection{Learning Latent Variables with Autoencoders}
Autoencoders (AEs) provide a straightforward means for manifold learning, since they can simultaneously learn a complex (nonlinear) mapping $c(\mathbf{x}^*)$ and an associated inverse mapping $c^{inv}(\bm{\phi})$ using simple basis functions. Linear autoencoders operate in the same space as PCA \cite{goodfellow2016}, while Kramer \cite{kramer1991} demonstrated the effectiveness of nonlinear autoencoders as a form of nonlinear PCA. A brief overview of the technique as relevant to the current work is presented here; the interested reader is referred to Chapter 14 of the book by Goodfellow et al. \cite{goodfellow2016} for a discussion of autoencoders, their uses, and comparisons to other manifold learning techniques.

Briefly speaking, an \textit{autoencoder} is a feed-forward artificial neural network that aims to replicate its input at its output. At a particular hidden layer within the autoencoder, the input is described as a ``code'', or $\bm{\phi}$. The dimensionality of the code $\bm{\phi}$ is determined by the structure of the neural network. The full network represents $(c^{inv} \circ c)(\mathbf{x}^*)$. The autoencoder is naturally split into the layers leading to $\bm{\phi}$, or $\bm{\phi} = c(\mathbf{x}^*)$, and the subsequent layers, or ${\mathbf{x}^*}' =  c^{inv}(\bm{\phi})$. The output of the network, ${\mathbf{x}^*}'=c^{inv}(c(\mathbf{x}^*))$ is an estimate of the original input $\mathbf{x}^*$. The autoencoder is typically trained using an iterative method by minimizing a loss function $L$,
\begin{eqnarray}
L(\mathbf{x}^*, {\mathbf{x}^*}'=c^{inv}(c(\mathbf{x}^*)))
\end{eqnarray}
which penalizes discrepancies between ${\mathbf{x}^*}'$ and $\mathbf{x}^*$. Commonly used loss functions include the mean squared error, mean absolute error, and variations of the hinge and cross-entropy functions.

\begin{figure}
	\centering
	\includegraphics[width=7cm]{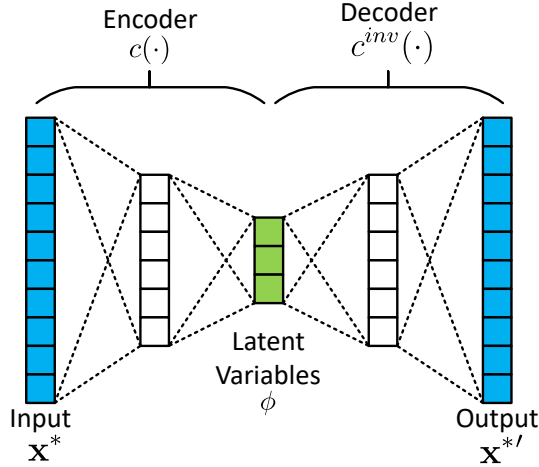}
	\caption{Conceptual depiction of an undercomplete autoencoder.}
	\label{fig:autoencoder}
\end{figure}

For the purpose of learning a low-dimensional manifold underlying a set of input data, we are particularly interested in \textit{undercomplete autoencoders}, or those with $\bm{\phi}$ constrained to have a lower dimension than $\mathbf{x}^*$. By restricting, or ``bottlenecking'', information flow through the feed-forward neural network, undercomplete autoencoders capture the salient trends present in the training data. Figure \ref{fig:autoencoder} depicts the structure of an undercomplete autoencoder with a two-layer encoder, two-layer decoder, and an encoded dimensionality of three. Undercomplete autoencoders are often constructed with an encoder and decoder that each comprise a single hidden layer. The universal approximator theorem \cite{hornik1991} guarantees that a feedforward neural network with at least one hidden layer can approximate any function (within a broad class) to an arbitrary degree of accuracy, given that enough hidden units are present. In practice, however, autoencoders with multiple hidden layers (termed \textit{deep} autoencoders) can sometimes reduce the computational cost of representing certain functions, improve data compression, and/or decrease the amount of training data required \cite{goodfellow2016, hinton2006}. Thus, with enough hidden units through depth or breadth, any (nonlinear) mapping between $\mathbf{y}=h(\mathbf{x})$ and $\mathbf{u} = K(\mathbf{y}_{sp} - h(\mathbf{x}))$ as relevant to the dynamical system under consideration here can be modeled to arbitrary accuracy, provided that $K(\cdot)$ is bounded and continuous. Note that MPC may not always satisfy this property, and alternative manifold learning techniques may be more suitable for systems exhibiting several distinct operating regimes.

\subsection{Building Latent Variable SBMs}

After a latent manifold underlying the closed-loop dynamics of a process is learned, the process operating data can be projected to the latent manifold to produce a low-dimensional representation. In particular, each observation $\mathbf{x}^*(t)$ can be transformed to $\bm{\phi}(t) = c(\mathbf{x}^*(t))$. Then, given a data set (e.g., the data set used for manifold learning) of transformed observations, $\bm{\phi}(t)$, and process setpoints, $\mathbf{y}_{sp}(t)$, we can perform system identification in the latent variable space to create a scale-bridging model (SBM) of the latent variable dynamics, $\dot{\bm{\phi}} = f^\phi(\bm{\phi},\mathbf{y}_{sp})$.
The system identification step can introduce additional inaccuracy in the dynamics embedded in the latent variable scheduling problem \eqref{eq:latentsched1}--\eqref{eq:latentsched2}; however, this is equally true when identifying SBMs using physical process variables, as in \eqref{eq:sbm1}--\eqref{eq:sbm2} \cite{pattison2016, tsay2019}. The identification of accurate dynamic models is a crucial step in both data-driven approaches and is performed using the same methods in either case.

Previous works \cite{pattison2016, tsay2019} employed SBMs in the Hammerstein-Wiener (HW) form to capture the closed-loop process dynamics of actual, physical process variables for scheduling applications. This choice was motivated by the inherent structure of HW models. In contrast to unstructured dynamic models (e.g., recurrent neural networks \cite{goodfellow2016}), HW models have fewer parameters and may be trained with significantly lower amounts of data. This is an important feature, since system identification experiments carried out on chemical plants can be expensive and time-consuming. A HW model comprises a linear dynamic component flanked by static, nonlinear input and output transformations. A single-input single-output (SISO) HW model can be written with the linear dynamic component represented as a state-space model:
\begin{align}
h &= H(y_{sp}) \label{eq:hw1} \\
\dot{\vec{r}} &= A\vec{r} + Bh \\
w &= C \vec{r} \\
x^* &= W(w) \label{eq:hw2}
\end{align}
where $H$ and $W$ are, respectively, the Hammerstein and Wiener blocks corresponding to the static, nonlinear input and output transformations. $A$, $B$, and $C$ are the matrices defining the linear state-space dynamical system, which is of order $n_{d}$, with $\vec{r} \in \R^{n_{d}}$. The SISO model in \eqref{eq:hw1}--\eqref{eq:hw2} is written for a single model input $y_{sp}$ and a single model output $x^*$.

The order of the linear dynamics and the choice of nonlinear transformations represent structural decisions that can also be made based on the data, such as by employing some information criterion \cite{tsay2019}. Typical nonlinear transformations for $H$ and $W$ include piecewise linear functions, sigmoid networks, saturation functions, or polynomials \cite{matlab}. Once the order of the linear state-space system and the nonlinear transformations are determined, the parameters of the $H$, $W$, and the dynamical system are fitted simultaneously to a training data set using an iterative algorithm. For systems with multiple production setpoints, a multi-input, single-output (MISO), HW model can be identified for each system output \cite{tsay2019}. Since HW models were successfully applied in previous works, we employ HW models to model the dynamics of latent variables in this work. The selection of HW models for this study also facilitates comparison of system identification and scheduling results to related works \cite{pattison2016, pattison2017, kelley2018, dias2018, tsay2018, tsay2019} involving data-driven models in the HW form.

\section{Case Study: Demand-Response Scheduling of an Air Separation Unit}
	
\subsection{Description of ASU Process}
Numerous previous works have considered the optimal scheduling of cryogenic air separation units (ASUs). Due to their large electricity consumption, ASUs can derive significant economic benefits from scheduling production in response to time-varying electricity prices (demand response). A common approach for scheduling ASUs is to assume quasi-stationary operation and use additional constraints to reflect the transition capabilities of the plant \cite{zhang2015, zhang2016b} (i.e., its ability to change its momentary power consumption by altering production rates/products); a few recent examples are reviewed here. Zhang et al.\cite{zhang2016b} employed a mode-based scheduling framework with surrogate sub-process models for computationally efficient scheduling calculations. Similarly, Zhou et al. \cite{zhou2017} defined a set of operational modes from historical data and used associated convex hulls for schedule optimization. More recently, Zhao et al. \cite{zhao2019} proposed a state-transition network model for scheduling ASUs, and applied it to a large-scale scheduling problem including two multi-product ASUs. Obermeier et al. \cite{obermeier2019} defined a mode-based scheduling approach to examine the important relationship between DR scheduling and equipment fatigue.

\begin{figure}
	\centering
	\includegraphics[width=10cm]{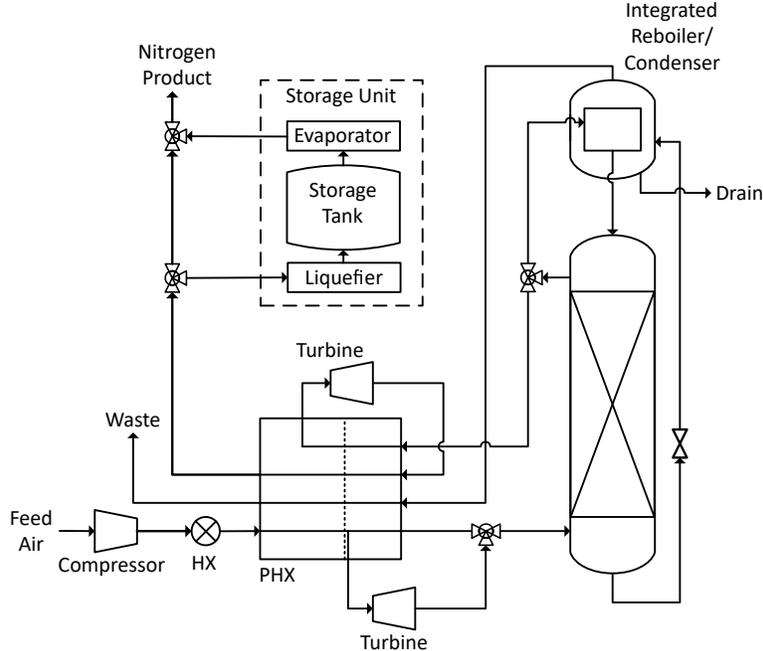}
	\caption{Flowsheet for a small nitrogen-production ASU with liquefier and liquid storage capacity.}
	\label{fig:ASUflowsheet}
\end{figure}

In this work, we consider the single-product ASU as shown in Figure \ref{fig:ASUflowsheet}. The detailed mathematical model of the process dynamics is based on the work of Cao et al. \cite{cao2015}, and is presented in full by Johansson \cite{johansson2015}. Pattison et al. \cite{pattison2016} investigated the scheduling problem using a full-order, detailed process model of the type \eqref{eq:scheduling1}--\eqref{eq:scheduling2}, as well as the SBM scheduling problem \eqref{eq:sbm1}--\eqref{eq:sbm2}. Dias et al. \cite{dias2018} developed a MPC system for the process and applied a novel simulation-optimization framework for integrated scheduling and control including MPC. Here, we employ the MPC system and its associated state-space models given by Dias et al. \cite{dias2018} with slight modification. The mathematical model of the process and its control system are summarized below.

The process in Figure \ref{fig:ASUflowsheet} produces high-purity nitrogen from an inlet air feed stream. The feed stream is compressed from atmospheric pressure to 6.8 bar, cooled, and passed through a primary (multi-stream) heat exchanger (PHX) where it is condensed against warming cryogenic streams. A portion of the air is removed from the PHX at an intermediate point and is sent to a turbine to generate electricity; the remainder exits the PHX at its saturation point. The two streams are combined and sent to the bottom of a cryogenic distillation column, which separates nitrogen from the other components of air. The bottoms product of the column is expanded through a valve before entering the reboiler. The reboiler and condenser are integrated in a single unit, allowing the bottoms stream to provide cooling duty to the condenser. The distillate of the column comprises the desired high-purity nitrogen stream. A portion of the distillate is sent to the condenser and becomes the column reflux, while the remaining product stream is expanded in a second turbine after being vaporized in the PHX. The product stream and the waste nitrogen stream from the reboiler both pass through the PHX to provide cooling duty to the incoming air. The nitrogen liquefier, storage, and evaporator units are included in the flowsheet. These units allow the plant to liquefy and store excess gaseous nitrogen generated during periods of over-production, and conversely evaporate stored liquid nitrogen to satisfy gas nitrogen demand during times of under-production.

The full-order process model comprises 6,094 equations and has 430 differential variables. The entire model is implemented in gPROMS, and implementation details can be found in our previous works \cite{dias2018, pattison2016}. The ASU process is assumed to operate with a constant gas nitrogen demand of 20 mol/s, with less than 2000 ppm impurity content (oxygen and argon). The process is assumed to be able to modulate its production rate by $\pm$20\% from its nominal value, representing a production range of 16 mol/s to 24 mol/s. The MPC for the process has four controlled variables and three manipulated variables. While the liquid drain stream from the reboiler was previously used as a fourth manipulated variable \cite{dias2018}, we found that outputs are not sensitive to this input in the desired range of operation. The MPC variables are summarized in Table \ref{tab:mpc}. The MPC has a sample time of six minutes and employs a linear state-space model created from system identification tests on the full-order dynamic model. The production rate setpoint represents $\mathbf{y}_{sp}(t)$, and its profile is set by the solution to the scheduling problem. The setpoints for the remaining controlled variables are fixed at $I_p^{sp}$ = 500 ppm, $\Delta T_{IRC}^{sp}$ 2.2 K, and $M_{reb}^{sp}$ = 100 kmol.

\begin{table}
	\centering
	\begin{tabular}{c c}
		Controlled Variable $y_p$ & Manipulated Variable $u_p$ \\ \hline
		Production flow rate & Inlet air flow rate \\
		Product impurity & PHX split fraction \\
		IRC temperature difference & Vapor product split \\
		Reboiler liquid level &
	\end{tabular}
\caption{Summary of MPC variables for the ASU Process}
\label{tab:mpc}
\end{table}

\subsection{Simulation Strategy for Generating Training and Testing Data}
\label{sec:datasim}
The augmented state variable vector for scheduling the ASU comprises 15 variables, i.e., $\text{dim}(\mathbf{x}^*) = 15$: the seven variables of the MPC (Table \ref{tab:mpc}), the power consumption, and seven state variables that feature constraints. These are the storage level M$_\text{store}$, column weeping ratio, column flooding ratio, column sump level, bubble-point pressure ratio, dew-point pressure ratio, nitrogen pressure ratio.

The column weeping ratio is defined as the minimal stage-wise ratio of vapor velocity to weeping velocity, while the flooding ratio is defined as the maximum stage-wise ratio of vapor velocity to flooding velocity. The bubble-point pressure ratio, defined as the ratio of pressure to bubble-point pressure for the stream exiting the PHX must be greater than one to ensure the stream is fully liquefied. The dew-point and nitrogen point pressure ratios, defined as the ratio of pressure to dew-point pressure for, respectively, the feed stream drawn at an intermediate location in the PHX and the product stream passing through the turbine, must be less than one to ensure the streams are fully vapor-phase. Note that $\mathbf{x}^*$ contains several each of input variables, state variables, and output variables. The manipulated variables are included in $\mathbf{x}^*$ to understand the degree to which modulating plant operations is possible and potential effects on the equipment \cite{obermeier2019, wiebe2018}. A full list of the variables included in $\mathbf{x}^*$, as well as their omission/inclusion status in previous studies \cite{pattison2016, dias2018}, is provided in the Appendix.

\begin{figure}[!h]
	\centering
	\includegraphics[width=16cm, trim={2cm, 0, 1cm, 0}, clip]{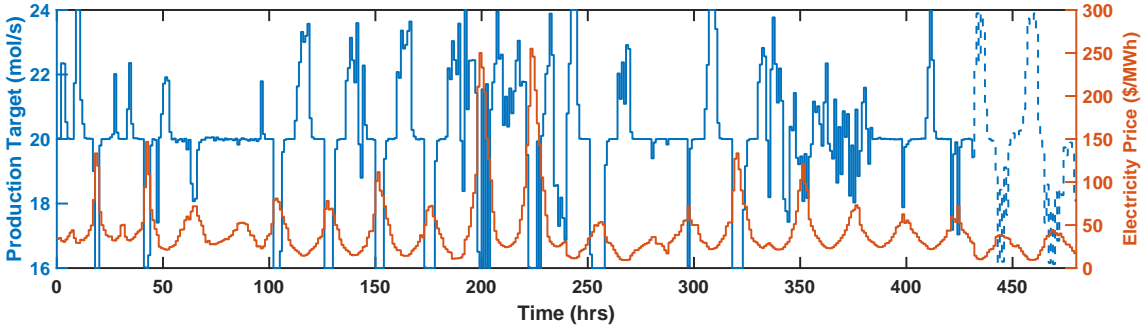}
	\caption{Electricity prices and production setpoints used to generate data set. The dashed lines depict the last 10\% of the data, which are used for dynamic model validation.}
	\label{fig:dataset}
\end{figure}

A data set was simulated for manifold learning using the detailed first-principles process model described by Johansson \cite{johansson2015}. The MPC was implemented ``online'' by linking the full-order dynamic model with the Model Predictive Control Toolbox in MATLAB. The full-order model was run between each MPC interval to generate sampled state variable values, and the MPC problem was solved in MATLAB to provide updated setpoints for the local regulatory controller in the subsequent interval. To generate an operating data set that reflects production modulation, the SBM-based scheduling problem \eqref{eq:sbm1}--\eqref{eq:sbm2} was solved using two-day electricity price data from a regional independent service operator (ISO). Ten such two-day price signals were selected, aiming to include a wide gamut of prices and hence process closed-loop behaviors. In total, 20 days of operating data were included in the data set. The electricity prices and resulting production targets used to generate the data set are shown in Figure \ref{fig:dataset}. Static dimensionality reduction techniques were applied to the simulated data set, and their statistics, as presented below, were computed using 5-fold cross-validation.

\begin{figure}[!h]
	\centering
	\includegraphics[width=8cm]{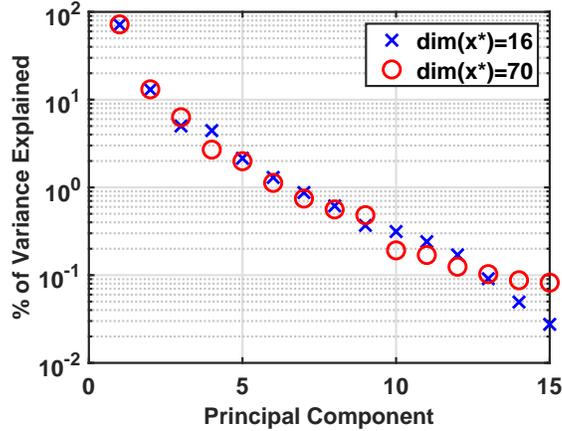}
	\caption{Percentage of variance in the full data set explained by each principal component.}
	\label{fig:pcaratios}
\end{figure}

Principal component analysis (PCA) reveals the linear relationships and correlations present in $\mathbf{x}^*$ in the data set. The percentage of variance explained by each principal component for the full data set is shown in Figure \ref{fig:pcaratios}, and the first few component loadings are shown in the Supporting Information. To ensure the correlation indicated by PCA is not coincidental, the same analysis was applied to the case where $\mathbf{x}^*$ included all \textit{process-level} variables. Here, the set of process-level variables refers to the properties (i.e., temperature, pressure, composition) of inter-unit streams and the operating conditions of the process units. Excluding variables that are necessarily identical to others and those with fixed/set values, there are 70 process-wide variables in total. Note that while the full dynamic model includes 6,094 equations, we limited this analysis to process-level variables, which provide ample information for most scheduling calculations. For example, many of the variables present in the model are associated with the spatial discretization of the primary heat exchanger and would not be measured in real-time in the plant.

PCA of the full data set produced a similar result to the case of 15 variables. In both cases, the percentage of variance explained decreases quickly with an increasing number of principal components (note the logarithmic ordinate scaling in Figure \ref{fig:pcaratios}), suggesting that the closed-loop process dynamics of $\mathbf{x}^*$ can be ``collapsed'' to a lower dimension, and that the accuracy of the approximation can be tuned by carefully selecting the dimensionality of the reduced-order representation. The variance captured decays more slowly after approximately 12 components for the case of all process-wide variables, since only the first 15 components (of 70) are shown.

\subsection{Manifold Learning Results}
Though PCA showed that the data set could be approximated reasonably using a low-dimensional set of linear latent variables, we incorporate nonlinear manifold learning techniques to further capture the closed-loop process dynamics. Manifold learning on the data set was performed using autoencoders (AEs). Several AE architectures were tested: \textit{Tanh(2x)}, having \texttt{tanh} activation functions and one hidden layer in the encoder and decoder; \textit{Tanh(3x)}, having \texttt{tanh} activation functions with two hidden layers in the encoder and decoder; and \textit{Linear}, with linear activation functions. As mentioned above, the representation power of an AE can be increased through increasing the depth or breadth of the neural network. In this study, the breadth of hidden layers is fixed to the truncated average of the input dimension and the encoded dimension. The effect of increasing the network size is investigated by switching from a single hidden layer to two hidden layers. Note that adding hidden layers to a neural network with only linear transformations does not increase the representation power of the model, since linear combinations of linear basis functions remain linear.

\begin{figure}[!h]
	\centering
	\includegraphics[width=8cm]{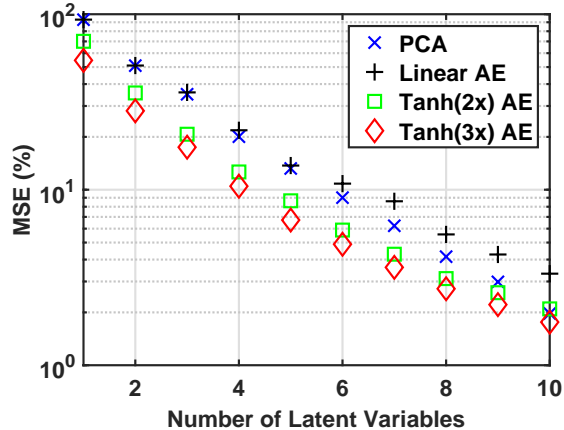}
	\caption{Comparison of validation MSE from dimensionality reduction techniques on process variables. The plotted MSEs were computed with 5-fold cross validation.}
	\label{fig:modelmse}
\end{figure}

Each process variable was scaled to take values between 0-100\%, with 0\% representing its minimum value in the data set and 100\% representing its maximum. The AEs were implemented and trained using TensorFlow \cite{tensorflow} with the Adam solver \cite{kingma2014} and the mean squared error (MSE) as the loss function:
\begin{equation}
\text{MSE} = \frac{||x_{ref}-x||^2}{N_s}
\end{equation}
where $N_s$ is the number of samples. Each AE was trained until the validation loss function remained the same or increased for several straight epochs. The cross-validated test MSEs from training the AEs are shown in Figure \ref{fig:modelmse}. While errors in the predictions for each process variable had the same weights in the loss function for this study, the errors of individual variables could be weighted differently to prioritize accuracy in certain process variables.

The results confirm the observation from the initial PCA, that the closed-loop process dynamics can be ``collapsed'' to a lower dimension. For example, with ten dimensions, the tested methods can represent the complete input data very well, with MSEs around 2\%. PCA presents an adequate manifold learning baseline, and linear AEs operate in the same space \cite{goodfellow2016}; however, we found that the accuracy of linear AEs to be lower than PCA when the number of latent variables included was large. This deviation can be attributed to difficulty in training a large AE using a stochastic optimization procedure, while PCA computes the optimal solution analytically. The nonlinear AEs achieve lower MSEs than the linear AE and PCA in all cases, with the benefits being significant especially at lower manifold dimensionality. The nonlinear AEs are capable of learning more complex relationships present in the data \cite{kramer1991}, and the increase from one hidden layer to two layers further increases representation power. The benefit of nonlinear mappings diminishes as the number of latent variables increases.

\subsubsection{Effect of Measurement Noise}
\begin{figure}[!h]
	\centering
	\includegraphics[width=8cm]{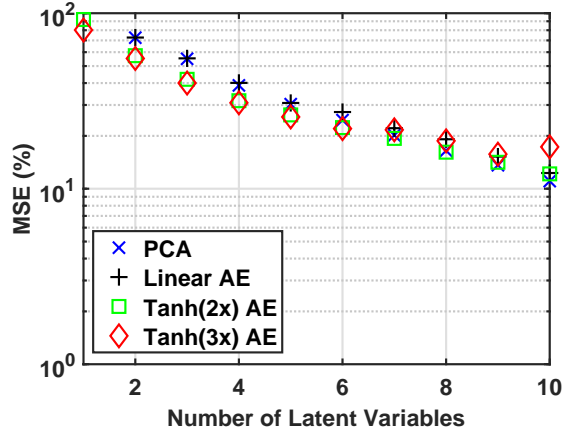}
	\caption{Comparison of validation MSE from dimensionality reduction techniques with 5\% normally distributed measurement noise. The plotted MSEs were computed with 5-fold cross validation.}
	\label{fig:noisymse}
\end{figure}

Dimensionality-reduction techniques are often employed for their ability to filter noisy data. Noise may sometimes be artificially introduced during autoencoder training to improve generalization ability \cite{goodfellow2016}. By constraining the intrinsic dimensionality of the retained information, latent variables retain the most important dynamics. In this case study, measurement noise was simulated by adding 5\% normally distributed error to all process variables in the training data set. AEs with the same configurations as described above were trained, and the MSEs are shown in Figure \ref{fig:noisymse}. Similar to the cases without measurement noise (Figure \ref{fig:modelmse}), the nonlinear AEs provide better process representations when the desired dimensionality is low, and the benefits of a nonlinear model decrease as more latent variables are added. All the models are less accurate when measurement noise is added, though this accuracy could be improved by increasing the size of the data set.

\begin{figure}[!h]
	\centering
	\includegraphics[width=8cm]{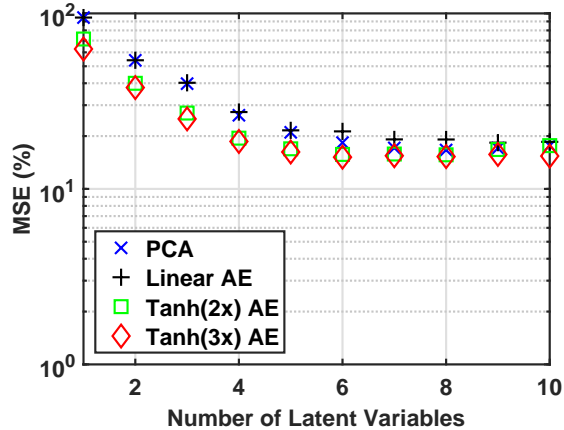}
	\caption{Comparison of validation MSE computed against ``ground truth'' from dimensionality reduction techniques with 5\% normally distributed measurement noise. The plotted MSEs were computed with 5-fold cross validation.}
	\label{fig:noisymse2}
\end{figure}

Insight into the denoising ability of the learned models can be obtained from examining their accuracy in predicting the original data set (i.e., the ``ground truth'' data without measurement noise). The MSEs of the models shown in Figure \ref{fig:noisymse} evaluated against the ground truth data are shown in Figure \ref{fig:noisymse2}. The MSEs decrease rapidly as the number of latent variables increases from one to six, where the model accuracy plateaus at approximately 11\% MSE. Interestingly, the learned models are generally more accurate in predicting the ground truth data than the noisy data (Figure \ref{fig:noisymse}), demonstrating their ability to filter noise. The introduction of measurement noise can be likened to a form of regularization, where shifting the variance-bias tradeoff improves model generalizability at the cost of some accuracy (the MSEs in Figure \ref{fig:noisymse2} are higher than those in Figure \ref{fig:modelmse}).  Since least-squares regression filters progressively more Gaussian noise as the number of samples increases, we again expect that the model accuracies could be improved by increasing the size of the data set.

\subsubsection{Effect of Additional Process Variables}
\label{sec:70vars}
\begin{figure}[!h]
	\centering
	\includegraphics[width=8cm]{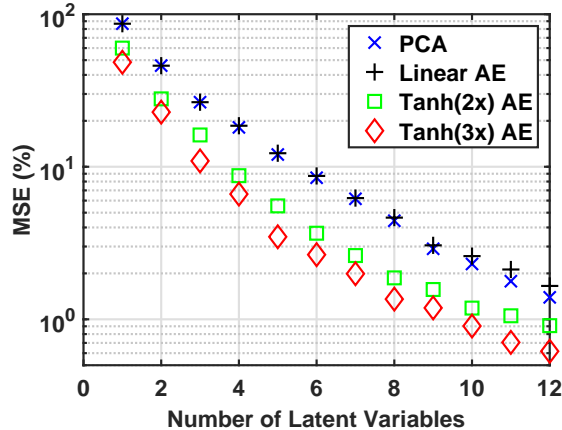}
	\caption{Comparison of validation MSE from dimensionality reduction techniques on all 70 process-level variables. The plotted MSEs were computed with 5-fold cross validation.}
	\label{fig:fullmse}
\end{figure}

To confirm that the low-order manifold mapping for $\mathbf{x}^*$ is not enabled by coincidentally selecting 15 correlated variables, the same AEs were trained on the full $\mathbf{x}^*$ vector that includes all 70 \textit{process-level} variables, as described in Section \ref{sec:datasim}. The cross-validated test MSEs from training the AEs on the full vector of process-wide variables are shown in Figure \ref{fig:fullmse}. The dimension of the hidden layers in the nonlinear AEs was again chosen to be the truncated average between the input layer and latent variable dimensions. The AEs were implemented and trained using the same procedure as above. Note that the AEs have more units (and representation power) due to the increase in dimension of the input layer.

As expected, increasing the dimension of $\mathbf{x}^*$ does not have a significant effect on the manifold learning procedure, since the intrinsic dimension $\text{dim}_i(\mathbf{x}^*)$ remains unchanged. This result supports the assertion that the dimensionality reduction is enabled by the low intrinsic dimension of the system. We again find that the closed-loop process dynamics of all 70 variables can be ``collapsed'' to a low dimension, with each additional dimension having a diminishing impact on model accuracy. The linear AE again exhibits a similar result to performing PCA on the data set, while the nonlinear AEs again perform better than both linear methods. The nonlinear AEs show improved accuracy compared to Figure \ref{fig:modelmse} since they have more hidden units. The MSEs for the nonlinear AEs reach $\sim$1\% with ten latent variables, while the linear models reach MSEs of around 2\%.

\subsection{Dynamic Modeling Results}
\begin{figure}[!h]
	\centering
	\includegraphics[width=16cm,trim={3cm, 0cm, 2cm, 0cm},clip]{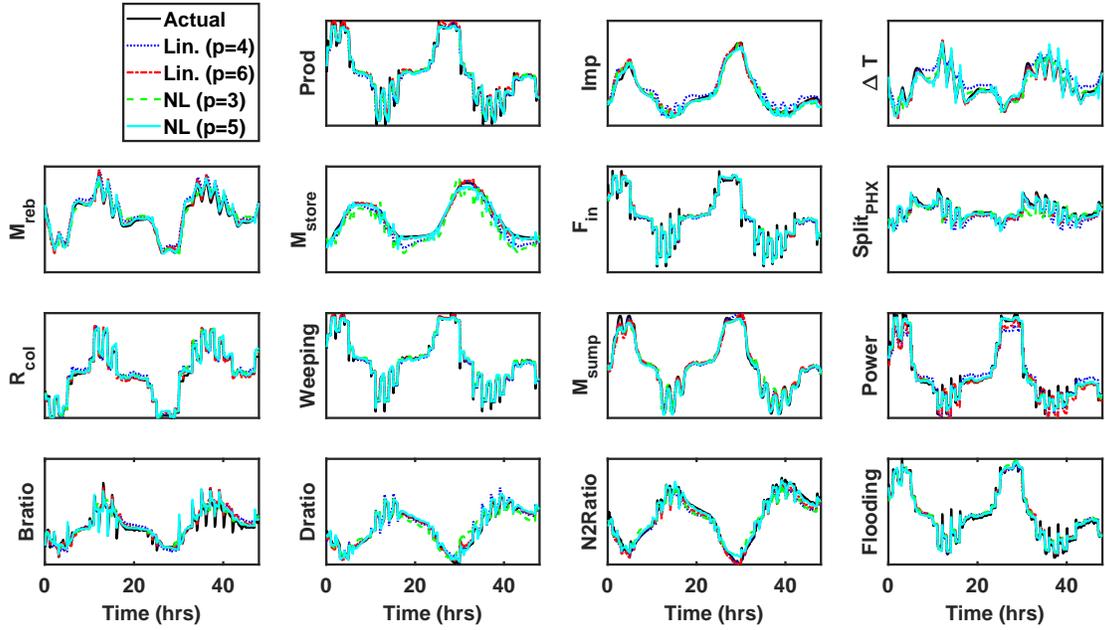}
	\caption{Evolution of the process variables predicted using various reduced-order representations, given ``true'' values of the latent variables.}
	\label{fig:testxprofiles_ae}
\end{figure}

\begin{figure}[!h]
	\centering
	\includegraphics[width=7cm]{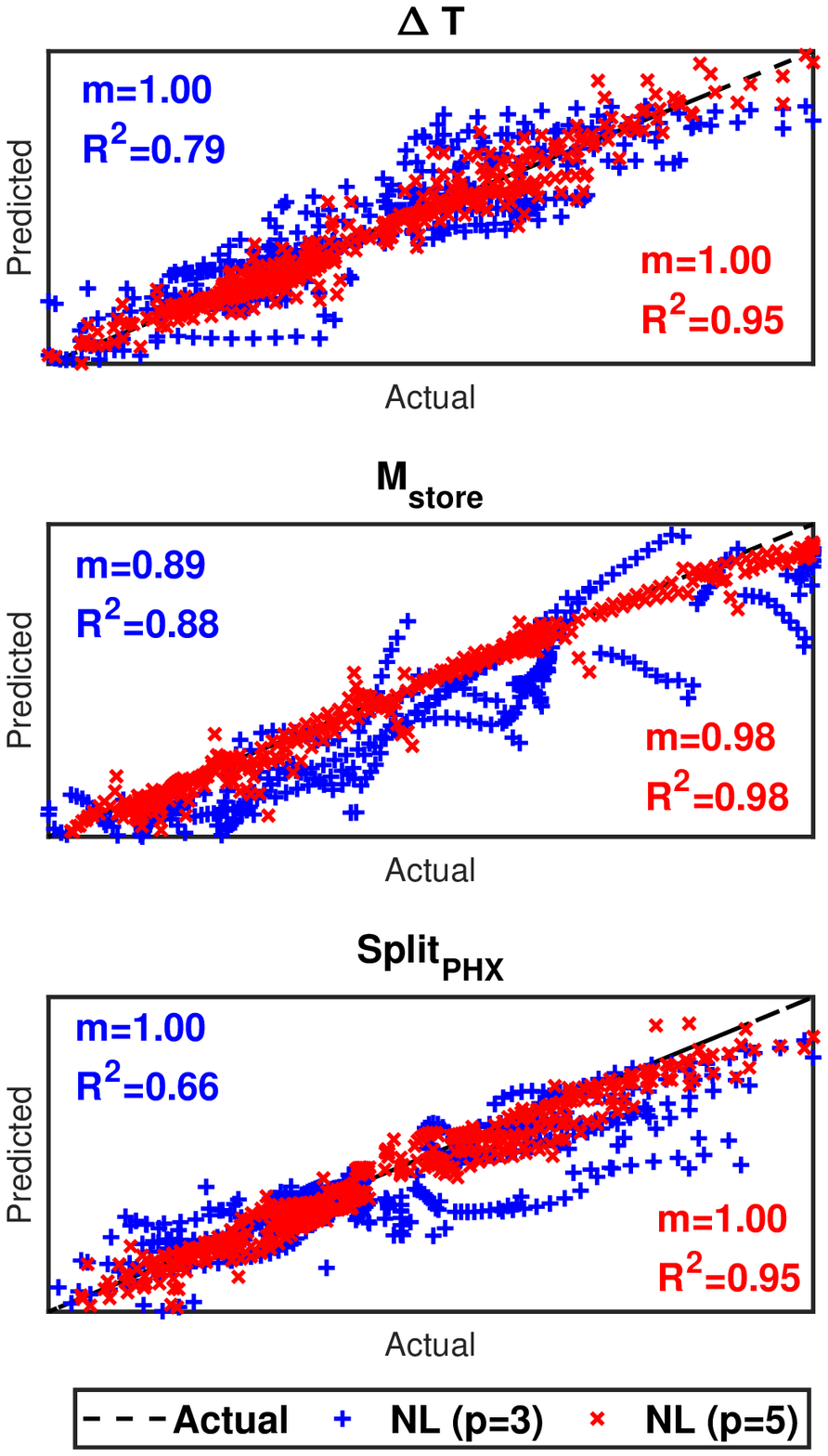}
	\caption{Scatter plots of a few process variable values predicted by nonlinear (NL) autoencoders with $p=3$ and $p=5$ latent variables. The predictions, ${\mathbf{x}^*}' = (c^{inv} \circ c)(\mathbf{x}^*)$, are plotted against their true values, $\mathbf{x}^*$, in the test data set.}
	\label{fig:scatter}
\end{figure}

Given the low MSEs possible with increasing $p$, we expect system identification to be the limiting factor in model accuracy for this study. Pattison et al. \cite{pattison2016} found that a 10\% ``back-off'' constraint was needed to compute feasible schedules for the ASU with HW models of physical variables, providing insight into SBM accuracy. We therefore select two low-dimensional representations of the ASU process dynamics with approximately 10\% MSE (Figure \ref{fig:modelmse}): a linear AE with six latent variables and a nonlinear AE with one hidden layer and five latent variables. To investigate the effect of adjusting $p$, we also test two representations with approximately 20\% MSE: a linear AE with four latent variables and a nonlinear AE with one hidden layer and three latent variables. For dynamic system identification, the first 18 days were used as training data (90\% of the data set). The process variables $\mathbf{x}^*$ were encoded using the respective encoders to give $\bm{\phi} = c(\mathbf{x}^*); \bm{\phi} \in \mathbb{R}^p$. The remaining two days are shown as dashed lines in Figure \ref{fig:dataset} and were left as test data. The test data were generated using electricity price data from a month not included in the training data to account for the potential for new patterns to emerge in the production schedules. The effect of dimensionality reduction on the test data was evaluated as a baseline by computing their estimated values ${\mathbf{x}^*}' = c^{inv}(\bm{\phi})$ using the true values of $\bm{\phi} = c(\mathbf{x}^*)$. The profiles of the variables $\mathbf{x}^*$ in the test data set, as well as their estimated values using all four AEs, are shown in Figure \ref{fig:testxprofiles_ae}. Scatter plots of ${\mathbf{x}^*}'$ against $\mathbf{x}^*$ are shown in Figure \ref{fig:scatter} for the most inaccurately predicted variables with the nonlinear AEs, where the improvement in prediction accuracy from increasing $p$ can easily be seen. Scatter plots for all variables in $\mathbf{x}^*$ are provided in the Supporting Information.

Hammerstein-Wiener (HW) models were fitted to the dynamics of the \textit{latent} variables $\dot{\phi}_i = f_i^\phi(\phi_i, \mathbf{y}_{sp})$. The models were fitted using the System Identification Toolbox in MATLAB \cite{matlab}. Piecewise-linear and polynomial transformations were used to model the $H$ and $W$ blocks of \eqref{eq:hw1}--\eqref{eq:hw2}. The form and polynomial order/number of piecewise segments for each transformation was determined by minimizing the normalized Akaike information criterion (nAIC) while using a large number of piecewise-linear segments for the other transformations and a high-order linear state-space model. The order of each linear state-space model was similarly determined using the nAIC. The resulting HW model structures and normalized mean squared error (NMSE) are shown in Tables \ref{tab:HWlinear} and \ref{tab:HWnonlinear}. Note that a higher NMSE indicates a better fit, in contrast to MSE. 
The NMSE values are computed using the \texttt{goodnessOfFit()} function in MATLAB, which uses the following formula:
\begin{equation}
\text{NMSE} = 1 - \frac{||x_{ref} - x||^2}{||x_{ref} - \text{mean}(x_{ref})||^2}
\end{equation}

\begin{table}[!h]
	\centering
	\begin{tabular}{c|c c c c|c c c c}
		& \multicolumn{4}{c}{Linear ($p=4$)} & \multicolumn{4}{|c}{Linear ($p=6$)} \\ \hline
		& $H$ &  & $W$ & NMSE & $H$ &  & $W$ & NMSE \\
		Variable & Form & $n_d$ & Form & train/test & Form & $n_d$ & Form & train/test \\ \hline
		$\phi_1$ & pwl-1 & 4 & poly-2 & 0.77/0.78 & pwl-3 & 5 & pwl-2 & 0.83/0.86 \\
		$\phi_2$ & pwl-2 & 4 & pwl-1 & 0.78/0.88 & pwl-1 & 6 & pwl-3 & 0.72/0.74 \\
		$\phi_3$ & pwl-2 & 5 & pwl-1 & 0.93/0.96 & pwl-1 & 4 & poly-2 & 0.90/0.93 \\
		$\phi_4$ & pwl-2 & 5 & poly-2 & 0.54/0.51 & pwl-4 & 4 & pwl-4 & 0.70/0.79 \\
		$\phi_5$ & - & - & - & - & poly-2 & 4 & pwl-3 & 0.85/0.90 \\
		$\phi_6$ & - & - & - & - & pwl-4 & 4 & pwl-5 & 0.64/0.63 \\ \hline
		average & - & - & - & 0.76/0.78 & - & - & - & 0.77/0.81
	\end{tabular}
	\caption{Hammerstein-Wiener model structures and accuracies for linear latent variables. Nonlinear transformations are denoted with `pwl' for piecewise-linear and `poly' for polynomial.}
	\label{tab:HWlinear}
\end{table}

\begin{table}[!h]
	\centering
	\begin{tabular}{c|c c c c|c c c c}
		& \multicolumn{4}{c}{Nonlinear ($p=3$)} & \multicolumn{4}{|c}{Nonlinear ($p=5$)} \\ \hline
		& $H$ &  & $W$ & NMSE & $H$ &  & $W$ & NMSE \\
		Variable & Form & $n_d$ & Form & train/test & Form & $n_d$ & Form & train/test \\ \hline
		$\phi_1$ & pwl-3 & 5 & pwl-3 & 0.76/0.83 & pwl-5 & 8 & pwl-2 & 0.49/0.38 \\
		$\phi_2$ & pwl-4 & 4 & pwl-2 & 0.77/0.89 & poly-2 & 8 & pwl-1 & 0.74/0.83 \\
		$\phi_3$ & pwl-3 & 6 & poly-3 & 0.53/0.70 & pwl-2 & 4 & poly-3 & 0.89/0.90 \\
		$\phi_4$ & - & - & - & - & pwl-3 & 4 & pwl-3 & 0.69/0.83 \\
		$\phi_5$ & - & - & - & - & pwl-2 & 5 & pwl-2 & 0.70/0.77 \\ \hline
		average & - & - & - & 0.69/0.81 & - & - & - & 0.70/0.74
	\end{tabular}
	\caption{Hammerstein-Wiener model structures and accuracies for nonlinear latent variables. Nonlinear transformations are denoted with `pwl' for piecewise-linear and `poly' for polynomial.}
	\label{tab:HWnonlinear}
\end{table}

\begin{figure}[!h]
	\centering
	\includegraphics[width=16cm,trim={3cm, 0cm, 2cm, 0cm},clip]{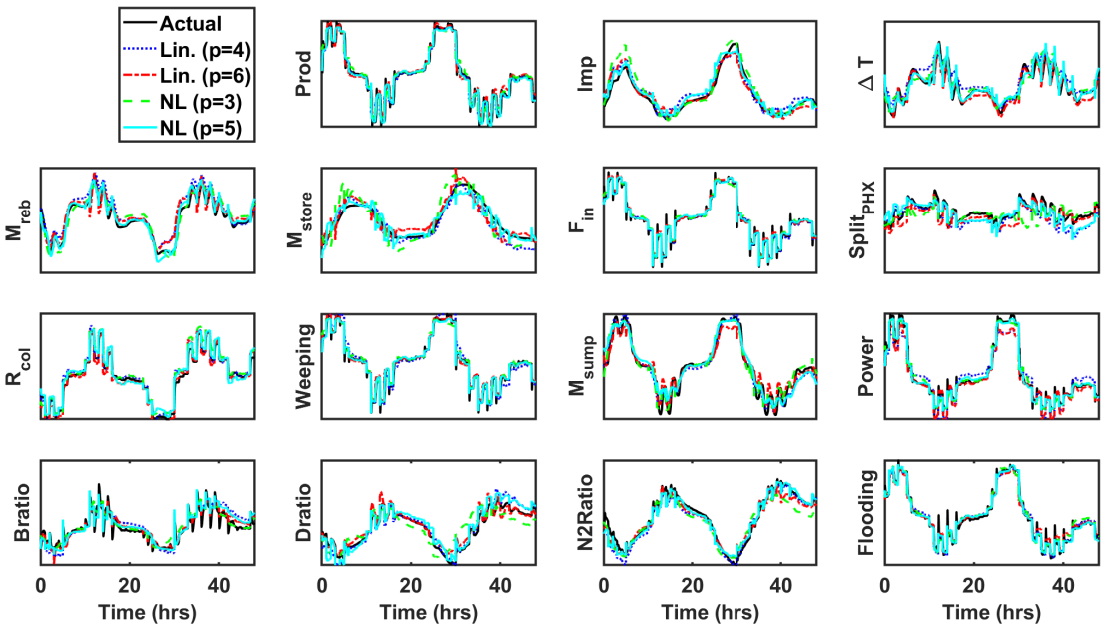}
	\caption{Evolution of the process variables predicted using various reduced-order representations, given values for the latent variables predicted by the identified HW models.}
	\label{fig:testxprofiles}
\end{figure}

The latent variable HW models were then simulated with the respective decoder $c^{inv}(\cdot)$ (created during autoencoder training) incorporated as additional static equalities. The simulations provide estimates of the latent variables $\bm{\phi}$ and decoded estimates of the process variables ${\mathbf{x}^*}' = c^{inv}(\bm{\phi})$. The actual variable profiles $\mathbf{x}^*$ and the estimated process variable profiles ${\mathbf{x}^*}'(t)$ are shown in Figure \ref{fig:testxprofiles}. The NMSEs for all 15 variables are shown in Table \ref{tab:testx}. While the SBMs for the latent variables exhibited lower NMSEs (Tables \ref{tab:HWlinear} and \ref{tab:HWnonlinear}), the final predictions for the process variables have NMSEs (Table \ref{tab:testx}) comparable to (or even higher than) previous works. Pattison et al. \cite{pattison2016} reported an average validation NMSE of 83.75\% using HW models to directly represent the behavior of eight physical variables.

The agreement between the actual values and the estimated profiles is generally very good, confirming the closed-loop process dynamics are approximated well by a data-driven model whose dynamics are confined to a low-dimensional, intrinsic manifold. As expected, increasing the dimensionality of the latent manifold from four to six in the linear case, and from three to five in the nonlinear case, improves the accuracy of the model predictions. The predictions of the integrated reboiler-condenser temperature difference $\Delta T_{IRC}$ and the PHX split fraction suffer from the largest inaccuracy. The predictions are slightly improved by increasing dimensionality (Figure \ref{fig:scatter}), but a comparison between Figures \ref{fig:testxprofiles_ae} and \ref{fig:testxprofiles} suggests that the error is primarily introduced by the system identification step. We note that the identified models are purely data-driven and may suffer from overfitting, particularly in cases where the training/test data set does not include the current operational situation. Implementation may benefit from a monitoring technique to detect whether the plant is entering an operating regime that has not been previously explored, at which point the data-driven models should be re-identified.

\begin{table}[!h]
	\centering
	\begin{tabular}{c|c c c c}
		& \multicolumn{4}{c}{NMSE} \\
		Variable & Lin ($p=4$) & Lin ($p=6$) & NL ($p=3$) & NL($p=5$) \\ \hline
		Production & 0.95 & 0.96 & 0.95 & 0.97 \\
		Impurity & 0.92 & 0.91 & 0.81 & 0.93 \\
		$\Delta$T$_{IRC}$ & \bf{0.63} & 0.90 & \bf{0.68} & 0.86 \\
		M$_\text{reb}$ & 0.83 & 0.89 & 0.86 & 0.94 \\
		M$_\text{store}$ & 0.86 & 0.85 & \bf{0.71} & 0.90 \\
		Air Flow & 0.98 & 0.98 & 0.98 & 0.97 \\
		PHX Split & \bf{0.21} & \bf{0.52} & \bf{0.10} & \bf{0.31} \\
		R$_\text{col}$ & 0.96 & 0.95 & 0.94 & 0.95 \\
		Weeping Ratio & 0.96 & 0.97 & 0.97 & 0.98 \\
		Sump Level & 0.95 & 0.94 & 0.93 & 0.93\\
		Power & 0.91 & 0.94 & 0.97 & 0.98  \\
		B Ratio & \bf{0.71} & 0.88 & 0.80 & 0.77\\
		D Ratio & 0.89 & 0.91 & \bf{0.70} & 0.86 \\
		N$_2$ Ratio & 0.93 & 0.95 & 0.87 & 0.96 \\
		Flooding Ratio & 0.98 & 0.97 & 0.97 & 0.98 \\ \hline
		Average & 0.85 & 0.90 & 0.83 & 0.89
	\end{tabular}
	\caption{NMSEs found on validation data set with linear (Lin) and nonlinear (NL) autoencoders at various levels of dimensionality reduction. Values below 0.75 are in bold.}
	\label{tab:testx}
\end{table}

\subsection{Optimal Scheduling Results}
To verify the application of the proposed framework to the ASE demand-response integrated scheduling and control problem, the latent-variable scheduling optimization problem \eqref{eq:latentsched1}--\eqref{eq:latentsched2} was solved for the demand response operational scenario considered by Dias et al. \cite{dias2018}. In this scenario, the storage tank is assumed to have a maximum capacity of 200 kmol of liquid nitrogen, with an initial inventory of 50 kmol. The inventory must be returned to at least its initial value at the end of the scheduling horizon to eliminate the potential of reporting false economic benefits derived from selling pre-existing inventory. The day-ahead electricity prices are assumed to be known over a 48-hour horizon, and there is a constant demand of 20 mol/s for the gas nitrogen product (equal to the nominal capacity of the plant). We consider here only the ``offline'' scheduling problem (with no re-scheduling); however, the proposed framework allows the scheduling problem to be solved quickly ($\sim$100-200s), which may benefit online scheduling techniques in the future \cite{pattison2017, beal2017}.

\begin{figure}[!h]
	\centering
	\includegraphics[width=16cm,trim={3cm, 0cm, 2cm, 0cm},clip]{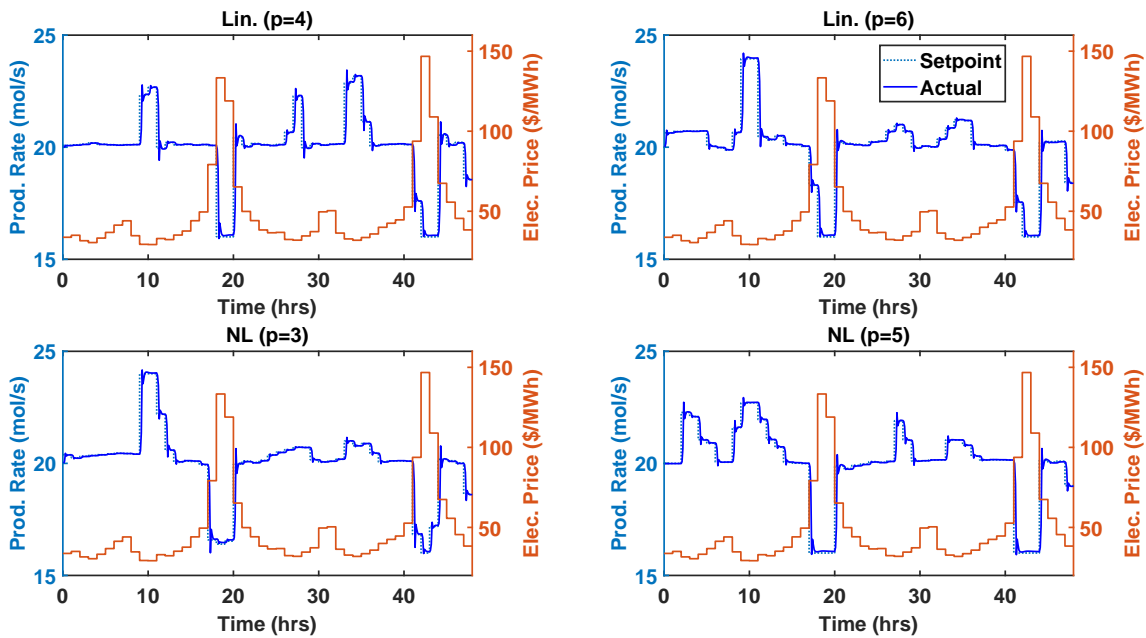}
	\caption{Electricity prices and the corresponding optimal production schedules.}
	\label{fig:schedules}
\end{figure}

\begin{table}
	\begin{center}
	\begin{tabular}{c c c c}
		& & Difference & Solution \\
		Case & Cost & from baseline & time \\ \hline
		Baseline$^a$ & \$707.91 & 0\% & - \\
		Simulation-optimization$^b$ & \$698.30 & 1.4\% & 381s \\
		Physical SBMs$^c$ & \$698.60 & 1.3\% & 610s \\
		Lin. ($p=4$) & \$700.75 & 1.0\% & 104s \\
		Lin. ($p=6$) & \$699.89 & 1.1\% & 113s \\
		NL ($p=3$) & \$701.65 & 0.9\% & 193s \\
		NL ($p=5$) & \$700.09 & 1.1\% & 238s \\
	\end{tabular}
	\end{center}
	\footnotesize
	$^a$ Baseline denotes the constant production rate case \\
	$^b$ Optimal point found using a simulation-optimization framework by \cite{dias2018}. Only variables involved in the MPC were modeled and constrained during optimization. \\
	$^c$ Optimal point found using SBMs identified for eight physical process variable, as reported by \cite{dias2018}. Details of the SBM models are provided by \cite{pattison2016}.
	\caption{Optimal schedule economic results}
	\label{tab:schedules}
\end{table}

\subsubsection{Linear Mappings}

The 48-hr scheduling problem was solved with the proposed latent-variable approach, using the linear mappings with $p=4$ and $p=6$. The models were implemented in gPROMS \cite{gproms}, and optimization was carried out using the built-in sequential dynamic optimization solver. The presented optimal points represent local optima found using 20 mol/s as the initial guess for the production setpoint at all times. The calculations were performed on a 64-bit Windows system with Intel Core i7-8700 CPU at 3.20 GHz and 16GB RAM. The implementation with four linear latent variables includes 29 differential variables and 64 total variables, while that with six linear latent variables includes 38 differential variables and 79 total variables. The scheduling problem with four latent variables was solved in 52 iterations, using 104.7s of CPU time (2.0s per iteration on average). The problem with six latent variables was solved in 42 iterations, using 113.0s of CPU time (2.7s per iteration on average). The two optimal schedules were simulated using the full-order dynamic model. The number of iterations taken by the local optimization solver to solve each problem is unpredictable. Nevertheless, the time per iteration is fairly consistent within each problem, and the number of iterations can be constrained for an expedited, but suboptimal, solution \cite{caspari2018}.

\begin{figure}[!h]
	\centering
	\includegraphics[width=16cm,trim={3cm, 0cm, 2cm, 0cm},clip]{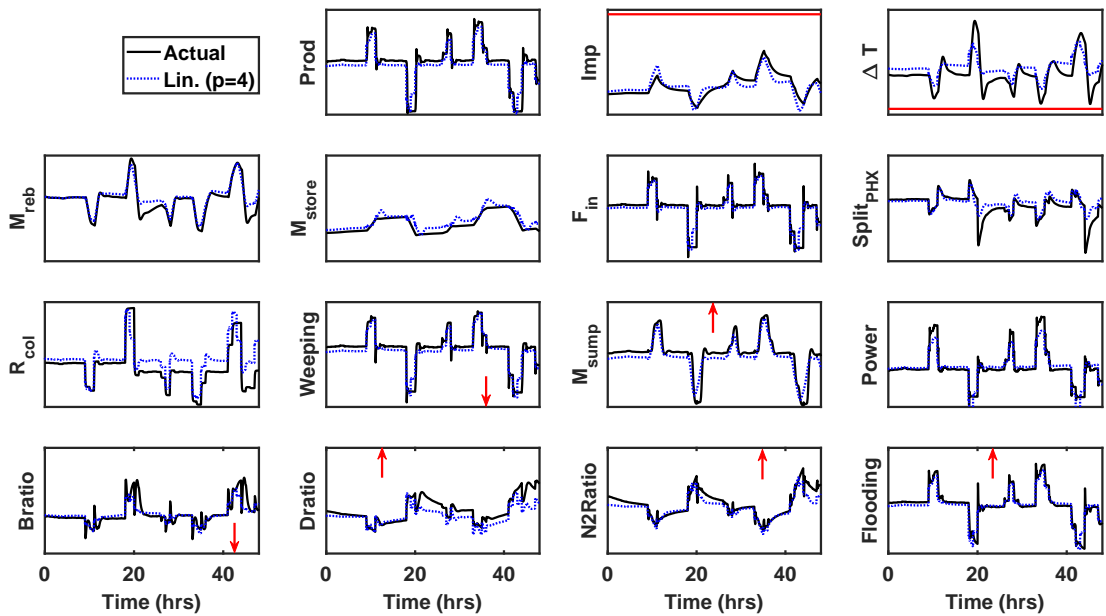}
	\caption{Optimal schedule generated with linear AEs (p=4). ``Actual'' profiles are generated by simulation of the same schedule using the full-order model with online MPC. Variable bounds are shown in red. Ordinate limits are 0-100\% of the respective scaled variable. Some bounds/constraints do not appear within this scaling; their locations are indicated with an arrow (e.g., an arrow pointing upwards indicates that an upper bound exists).}
	\label{fig:lin4}
\end{figure}

\begin{figure}[!h]
	\centering
	\includegraphics[width=16cm,trim={3cm, 0cm, 2cm, 0cm},clip]{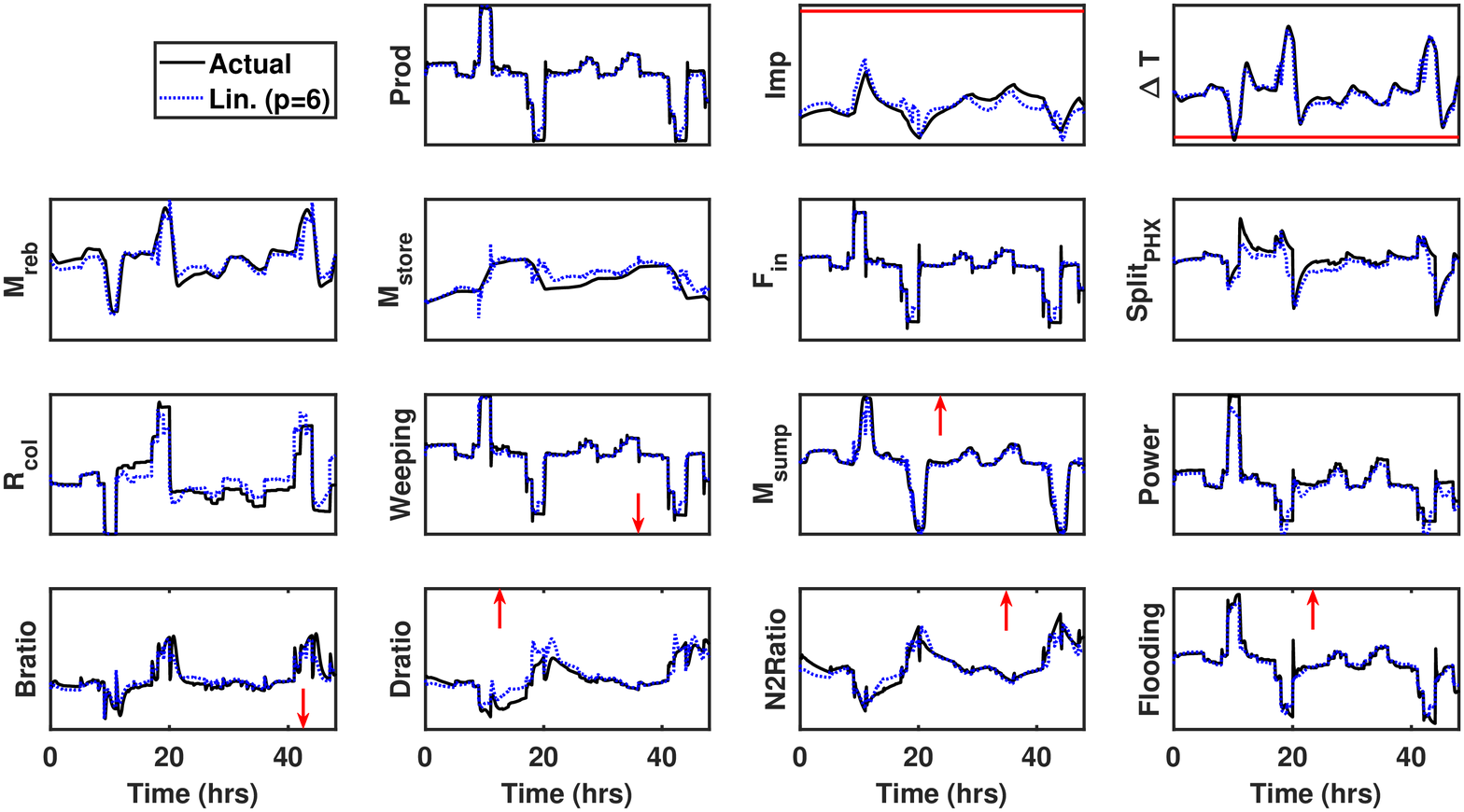}
	\caption{Optimal schedule generated with linear AEs (p=6). ``Actual'' profiles are generated by simulation of the same schedule using the full-order model with online MPC. Variable bounds are shown in red. Ordinate limits are 0-100\% of the respective scaled variable. Some bounds/constraints do not appear within this scaling; their locations are indicated with an arrow.}
	\label{fig:lin6}
\end{figure}

The production rates found in simulation of the optimal schedules are shown in Figure \ref{fig:schedules}, along with their setpoints/targets, which are closely tracked. As expected, production rates are scheduled to decrease when energy prices are high in both schedules. The behaviors of all 15 modeled process variables in the two computed schedules are shown in Figures \ref{fig:lin4} and \ref{fig:lin6}. The temperature driving force across the reboiler/condenser nearly reaches its bound in both schedules, but this potential constraint violation is only predicted by the low-dimensional representation when six linear latent variables are included. As shown in Table \ref{tab:testx}, increasing $p$ from four to six improves the test NMSE on $\Delta T_{IRC}$ from 0.63 to 0.90. None of the other variable path constraints were reached when the optimal schedules were simulated with the full-order dynamic model. The end point constraints of returning the storage and reboiler levels to at least their initial conditions were also met in both schedules. We note that although the variables with inactive constraints may not have been necessary for computing a feasible schedule, the proposed approach captures the dynamics of all constrained variables in the scheduling problem. This eliminates the difficult task of anticipating which constraints may be violated and should therefore be modeled.

The predictions are generally improved by including six (vs. four) latent variables, especially in the aforementioned temperature driving force $\Delta T$ and PHX split fraction. However, some deviations are visible for the column reflux $R_{col}$ predictions with both models (Figures \ref{fig:lin4} and \ref{fig:lin6}). The operational costs calculated using the full-order dynamic model of the computed schedules are shown in Table \ref{tab:schedules}. Both schedules result in a approximately 1\% savings compared to a constant production profile set at the nominal rate (subject to the same electricity price profile). These savings are similar to those reported in \cite{dias2018} and represent a substantial amount in the context of the well-established, commoditized air separation industry. In contrast to the previous approaches, the proposed method maintains the computational efficiency of scheduling in a reduced dimension while providing predictions of all constrained variables.

\subsubsection{Nonlinear Mappings}

The same 48-hr latent-variable scheduling problem was solved using the nonlinear mappings with $p=3$ and $p=5$. The same implementation and optimization settings were used, but the models include a hidden layer and nonlinear transformations.  The implementation with three nonlinear latent variables includes 25 differential variables and 66 total variables, while that with five nonlinear latent variables includes 38 differential variables and 86 total variables. The scheduling problem with three nonlinear latent variables was solved in 69 iterations, using 192.8s of CPU time (2.8s per iteration on average). The problem with five latent variables was solved in 61 iterations, using 237.8s of CPU time (3.9s per iteration on average). Although the problems with four linear and three nonlinear latent variables have a similar number of variables, the optimization problem with nonlinear latent variables requires more time per optimization iteration. The same phenomenon is observed for the problems with six linear and five nonlinear latent variables. The slowdown can be attributed to the nonlinear \texttt{tanh} transformations limiting integration step sizes and thereby slowing down the implicit time-integration scheme at each iteration.

\begin{figure}[!h]
	\centering
	\includegraphics[width=16cm,trim={3cm, 0cm, 2cm, 0cm},clip]{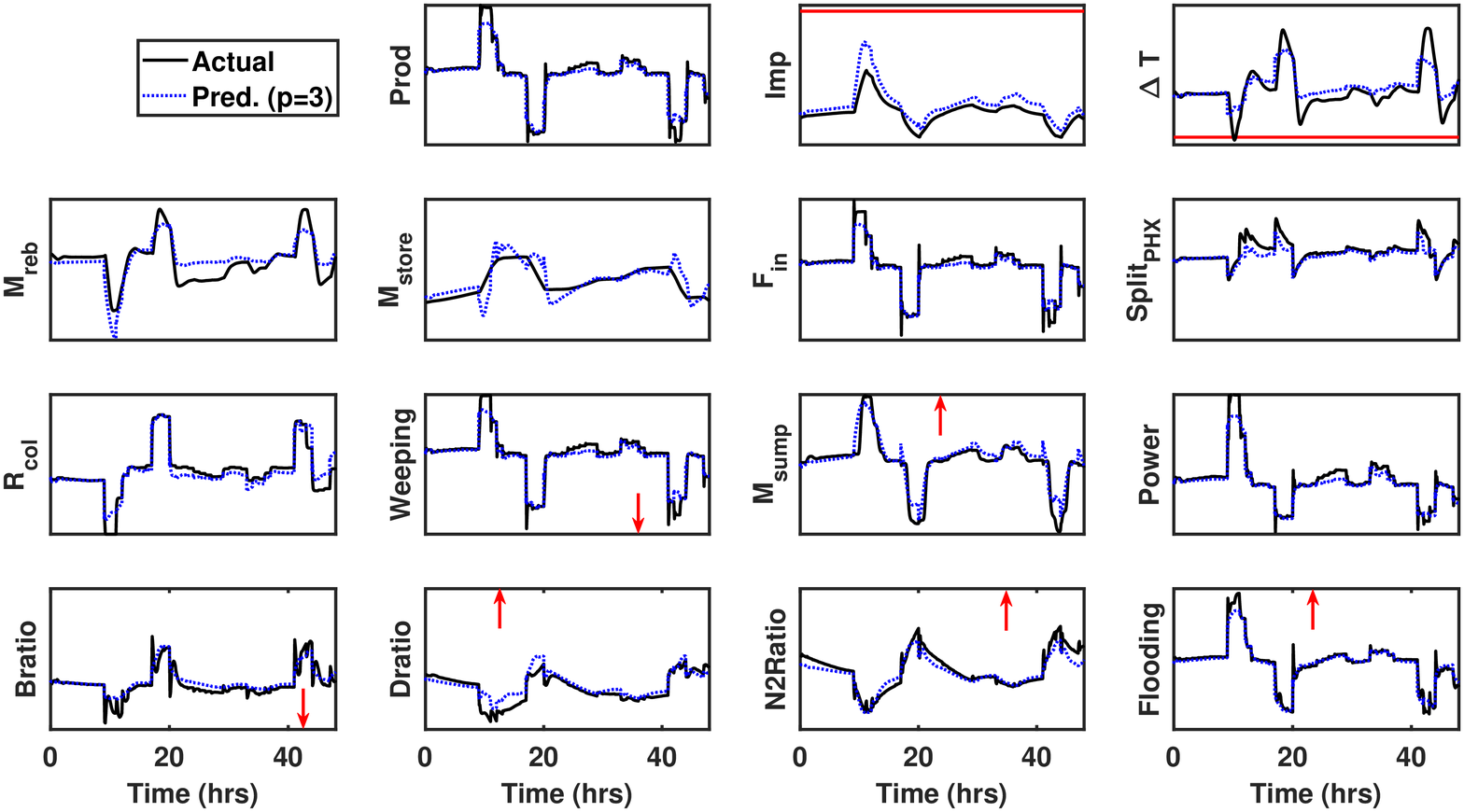}
	\caption{Optimal schedule generated with nonlinear AEs (p=3). ``Actual'' profiles are generated by simulation of the same schedule using the full-order model with online MPC. Variable bounds are shown in red. Ordinate limits are 0-100\% of the respective scaled variable. Some bounds/constraints do not appear within this scaling; their locations are indicated with an arrow.}
	\label{fig:nl3}
\end{figure}

\begin{figure}[!h]
	\centering
	\includegraphics[width=16cm,trim={3cm, 0cm, 2cm, 0cm},clip]{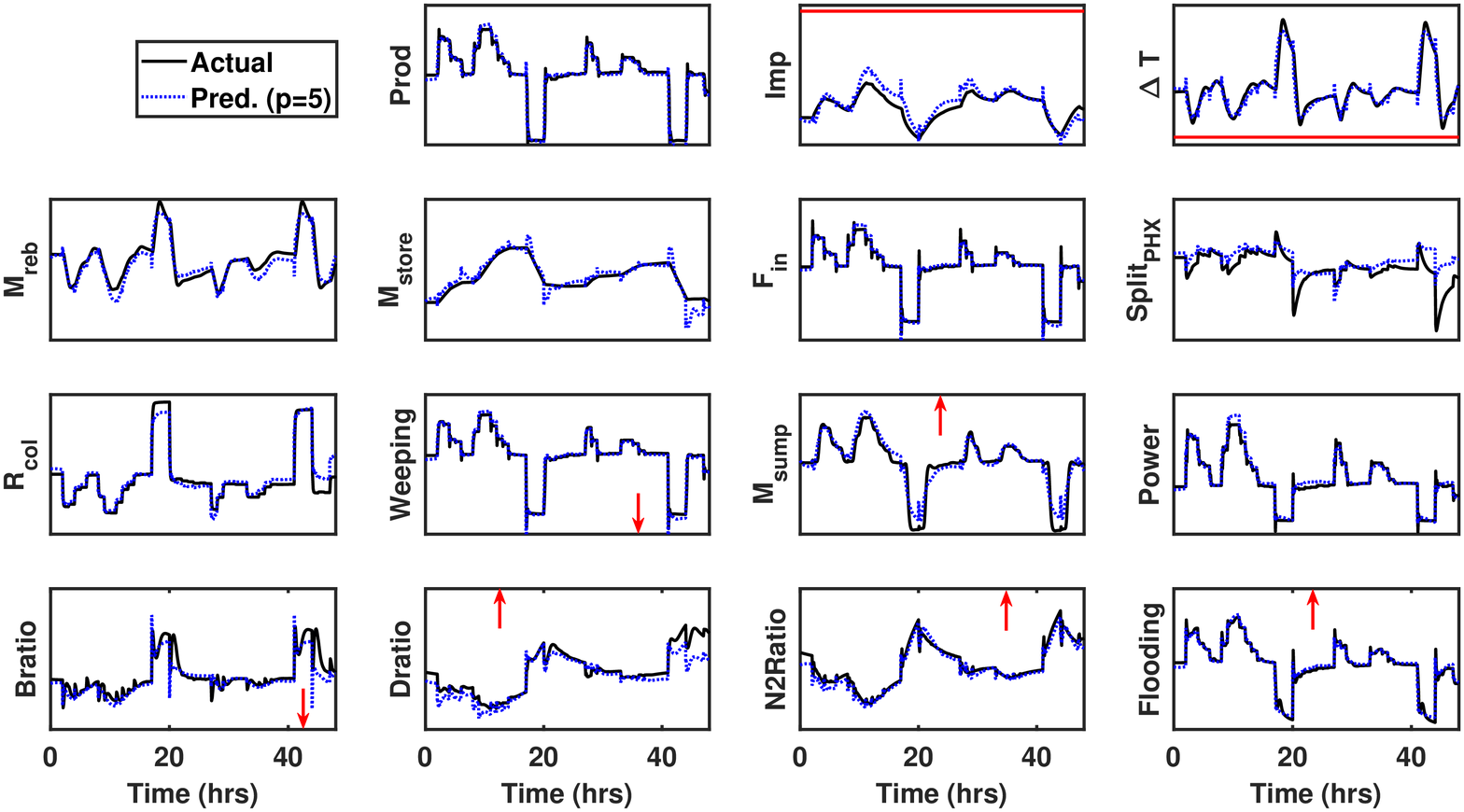}
	\caption{Optimal schedule generated with nonlinear AEs (p=5). ``Actual'' profiles are generated by simulation of the same schedule using the full-order model with online MPC. Variable bounds are shown in red. Ordinate limits are 0-100\% of the respective scaled variable. Some bounds/constraints do not appear within this scaling; their locations are indicated with an arrow.}
	\label{fig:nl5}
\end{figure}

The two computed optimal schedules were simulated with the aforementioned MPC and the full-order dynamic model. The production rate setpoints and actual values found at the optimal points are again shown in Figure \ref{fig:schedules}. The behavior of all 15 modeled process variables in the two computed schedules is shown in Figures \ref{fig:nl3} and \ref{fig:nl5}. The temperature driving force across the reboiler/condenser slightly violates the respective bound in the schedule computed with three nonlinear latent variables, which is not predicted accurately by the reduced-order model. Increasing $p$ from three to five improves the test NMSE on $\Delta T_{IRC}$ from 0.68 to 0.86 (Table \ref{tab:testx}), and the constraint violation is avoided by using five nonlinear latent variables. The end point constraints of returning the storage and reboiler levels to at least their initial conditions were met in both schedules.

The optimization problem with dynamics represented by five nonlinear latent variables demonstrates the most accurate predictions of the evolution of process variables, but also required the most time per optimization iteration out of the four tested latent-variable scheduling problems. However, the optimization problem size is still greatly reduced from previous approaches \cite{pattison2016, pattison2017, dias2018}, and correspondingly, the optimal schedule was obtained with less computational effort. In addition, the proposed formulation allows for more information on the process dynamics to be captured in scheduling calculations, with all constrained process variables predicted relatively accurately. The results in Section \ref{sec:70vars} further suggest that more process-level variables could be easily included at similar levels of accuracy without increasing the size of the latent dynamics. The intrinsic, low-dimensional dynamics underlying the closed-loop system would be approximated, and the decoder could be expanded to include more process variables.

\section{Conclusions}
The integration of closed-loop process dynamics in production scheduling calculations is key to ensuring that production schedules do not violate process constraints when implemented in practice (i.e., the schedules are ``dynamically feasible''). This is especially important when production changes are frequent, such as in practical cases driven by fast-changing external conditions. However, there is an intrinsic tradeoff between the amount of dynamic information captured in the scheduling model and the computational complexity required for its optimization. Driven by the need for computationally efficient representations of process dynamics, we exploit the low intrinsic dimensionality of closed-loop process behavior to generate reduced-order models. We proposed a data-driven approach for learning a low-dimensional latent manifold underlying variations in recorded observations, which can then be used to represent the process behavior. We present a conceptual analysis of the existence of such a manifold, and we demonstrate the means for selecting the dimensionality of the latent manifold, so as to balance between complexity of the captured dynamics and model size.

We presented a framework for production scheduling using the latent variable representation of process dynamics. In the proposed framework, process operating data are projected onto the latent manifold, and system identification and scheduling calculations are both performed in the latent variable space. Using this method, system identification is only necessary for a lower number of (latent) variables, and the scheduling calculations require less computational expense. These combined advantages allow for a broad complement of process variables to be modeled efficiently, and for the integrated scheduling and control problem to be solved in practically-relevant amounts of time. Notably, this approach eliminates the need for heuristically selecting variables to be modeled in scheduling calculations by automating the dimensionality reduction step. The framework was applied to the scheduling of an air separation unit under MPC in response to an hourly electricity price signal. The results confirm that the latent-variable scheduling formulation can retain more information about the process dynamics, compared to previous works, and simultaneously reduce the required computational effort.

\section{Disclaimer}
This report was prepared as an account of work sponsored by an agency of the United States Government. Neither the United States Government nor any agency thereof, nor any of their employees, makes any warranty, express or implied, or assumes any legal liability or responsibility for the accuracy, completeness, or usefulness of any information, apparatus, product, or process disclosed, or represents that its use would not infringe privately owned rights. Reference herein to any specific commercial product, process, or service by trade name, trademark, manufacturer, or otherwise does not necessarily constitute or imply its endorsement, recommendation, or favoring by the United States Government or any agency thereof. The views and opinions of authors expressed herein do not necessarily state or reflect those of the United States Government or any agency thereof.

\section{Acknowledgements}
The authors acknowledge funding from the National Science Foundation (NSF) through CAREER Award 1454433 and Award CBET-1512379. This material is also based on work supported by the US Department of Energy under Award Number DE-OE0000841. CT thanks The University of Texas at Austin for support through the University Graduate Continuing Fellowship.

\newpage
\noindent{\bf \textsc{Literature Cited}}

\small
\def\refname{}
\def\bibsection{}
\vspace{2mm}

\bibliographystyle{unsrt}
\bibliography{latentvarscheduling_CT3}

\newpage
\section*{Appendix}
\setcounter{table}{0}
\renewcommand{\thetable}{A\arabic{table}}
\setcounter{figure}{0}
\renewcommand{\thefigure}{A\arabic{figure}}
\vspace{-6pt}

\begin{table}[!h]
	\begin{center}
		\small
		\begin{tabular}{c c c c c c c}
			\# & Variable & Description & Unit & This Work & Pattison et al. \cite{pattison2016} & Dias et al. \cite{dias2018} \\ \hline
			1 & Prod & Production flow rate & mol/s & x & x & x \\
			2 & Imp & Product impurity & ppm & x & x & x \\
			3 & $\Delta$ T & Temperature difference in IRC & K & x & x & x \\
			4 & M$_\text{reb}$ & Reboiler liquid level & mol & x & x & x \\
			5 & M$_\text{store}$ & Storage liquid level & mol & x & x$^a$ & x \\
			6 & F$_\text{in}$ & Air flow in & mol/s & x & x & x \\
			7 & Split$_\text{PHX}$ & PHX split fraction & \% & x & - & x \\
			8 & R$_\text{col}$ & Vapor product split & \% & x & - & x \\
			9 & Weeping & Column weeping ratio & \% & x & x & - \\
			10 & M$_\text{sump}$ & Sump liquid level & mol & x & - & -\\
			11 & Power & Power consumption & MW & x & x$^a$ & x \\
			12 & Bratio & PHX outlet P / bubble point P & \% & x & x & - \\
			13 & Dratio & Turbine 1 inlet P / dew point P & \% & x & x & - \\
			14 & N2ratio & Turbine 2 inlet P / dew point P & \% & x & - & - \\
			15 & Flooding & Column flooding ratio & \% & x & - & -
		\end{tabular}
	\end{center}
	\footnotesize
	$^a$Pattison et al. \cite{pattison2016} leveraged process knowledge to estimate static models for the plant power consumption and storage level based on the production flow rate rather than create separate dynamic models.
	\caption{Variables dynamically modeled in ASU case study for various approaches.}
\end{table}

\begin{table}[!h]
	\begin{center}
	\small
	\begin{tabular}{c c c c c c c}
		\# & Variable & unit & Lower Bound & Upper Bound & 0\% Scaled Value & 100\% Scaled Value \\ \hline
		1 & Prod & mol/s & - & - & 15.81 & 24.17 \\
		2 & Imp & ppm & - & 1800 &158.68 & 1872.73 \\
		3 & $\Delta$ T & K & 1.9 & - & 1.85 & 2.80 \\
		4 & M$_\text{reb}$ & mol & 0 & - & 86148.62 & 109259.89 \\
		5 & M$_\text{store}$ & mol & 0 & - & 0.00 & 161097.47 \\
		6 & F$_\text{in}$ & mol/s & - & - & 30.41 & 49.62 \\
		7 & Split$_\text{PHX}$ & \% & - & - & 3.19 & 3.83 \\
		8 & R$_\text{col}$ & \% & - & - & 51.57 & 52.73 \\
		9 & Weeping & \% & 105 & - & 1902.67 & 3025.24 \\
		10 & M$_\text{sump}$ & mol & 0 & 5000 & 2002.46 & 2876.60 \\
		11 & Power & MW & - & - & 0.22 & 0.45 \\
		12 & Bratio & \% & 105 & - & 163.13 & 183.22 \\
		13 & Dratio & \% & - & 95 & 81.69 & 90.80 \\
		14 & N2ratio & \% & - & 95 & 57.85 & 69.01 \\
		15 & Flooding & \% & - & 95 & 68.68 & 96.57
	\end{tabular}
	\end{center}
	\caption{Scaled variable ranges and bounds for ASU case study. ``Back-off'' bounds are used to account for model inaccuracy.}
\end{table}

\end{document}